\documentclass{article}

\usepackage[utf8]{inputenc}
\usepackage[T1]{fontenc}
\usepackage[french,english]{babel}
\usepackage{amsmath}
\usepackage{amssymb}
\usepackage{amsthm}
\usepackage{amsfonts}
\usepackage{enumitem}
\usepackage{mathtools}
\usepackage{faktor}
\usepackage{geometry}
\usepackage{xparse}
\usepackage{tikz}
\usepackage{tikz-cd}
\usepackage{mathrsfs}
\usepackage{csquotes}
\usepackage{verbatim}
\usepackage{setspace} 
\usepackage[style=alphabetic,doi=true,backend=biber]{biblatex}
\DeclareFieldFormat{postnote}{#1} 
\usepackage{adjustbox}
\usepackage{hyperref} 
\usepackage{cleveref} 

\newtheorem{introthm}{Theorem}

\newtheorem{introcor}[introthm]{Corollary}

\newtheorem{thm}{Theorem}[section]
\newtheorem*{ass*}{Assumption}
\newtheorem{lem}[thm]{Lemma}

\newtheorem{cor}[thm]{Corollary}
\newtheorem{prop}[thm]{Proposition}
\newtheorem*{prop*}{Proposition}
\newtheorem*{thm*}{Theorem}
\newtheorem{defi}[thm]{Definition}
\newtheorem*{defi*}{Definition}

\theoremstyle{definition}

\theoremstyle{remark}
\newtheorem*{rmk}{Remark}

\numberwithin{equation}{section}

\crefname{thm}{thm.}{thms.}
\Crefname{thm}{Theorem}{Theorems}
\crefname{lem}{lem.}{lemmas}
\Crefname{lem}{Lemma}{Lemmas}
\crefname{prop}{prop.}{props.}
\Crefname{prop}{Proposition}{Propositions}
\crefname{defi}{def.}{defs.}
\Crefname{defi}{Definition}{Definitions}
\crefname{section}{§}{§}
\Crefname{section}{Section}{Sections}
\crefname{subsection}{§}{§}
\Crefname{subsection}{Subsection}{Subsections}
\crefname{paragraph}{§}{§}
\Crefname{paragraph}{Paragraph}{Paragraphs}
\crefname{cor}{cor.}{cors.}
\Crefname{cor}{Corollary}{Corollaries}

\newcommand*{\sheafhom}{\mathscr{H}\kern -1.5pt om}
\newcommand{\set}[1]{\left\{#1\right\}}
\newcommand{\ct}[1]{Z \kern -1 pt \left({#1}\right) \kern -0.5pt}
\renewcommand{\cal}{\mathcal}
\renewcommand{\ker}{\mathrm{ker}}

\newcommand{\ZZ}{\mathbb{Z}}
\newcommand{\NN}{\mathbb{N}}

\newcommand{\FF}{\mathbb{F}}
\newcommand{\QQ}{\mathbb{Q}}

\newcommand{\CC}{\mathbb{C}}
\DeclareMathOperator{\Spec}{Spec}
\DeclareMathOperator{\Hom}{Hom}

\DeclareMathOperator{\Ext}{Ext}

\DeclareMathOperator{\fib}{fib}

\DeclareMathOperator{\Gal}{Gal}

\DeclareMathOperator{\rank}{rank}

\DeclareMathOperator{\ord}{ord}
\DeclareMathOperator{\colim}{colim}

\DeclareMathOperator{\Perf}{Perf}

\DeclareMathOperator{\Cons}{Cons}

\tikzset{%
    symbol/.style={%
        draw=none,
        every to/.append style={%
            edge node={node [sloped, allow upside down, auto=false]{$#1$}}}
    }
} 

\newsavebox{\pullback}
\sbox\pullback{%
	\begin{tikzpicture}%
		\draw (0,0) -- (1ex,0ex);%
		\draw (1ex,0ex) -- (1ex,1ex);%
\end{tikzpicture}}

\title{Weil-étale cohomology and the equivariant Tamagawa number conjecture for constructible sheaves in characteristic \texorpdfstring{$p$}{p}}
\author{Adrien Morin}
\date{}

\begin{document}
	
\maketitle

\begin{abstract}

Let $X$ be a variety over a finite field. Given an order $R$ in a semi-simple algebra over the rationals and a constructible étale sheaf $F$ of $R$-modules over $X$, one can consider a natural non-commutative $L$-function associated with $F$. We  prove a special value formula at negative integers for this $L$-function, expressed in terms of Weil-étale cohomology; this is a geometric analogue of, and implies, the equivariant Tamagawa number conjecture for an Artin motive and its negative twists over a global function field. It also generalizes the results of Lichtenbaum and Geisser on special values at negative integers for zeta functions of varieties, and the work of Burns--Kakde in the case of non-commutative L-functions coming from a Galois cover of varieties.

\end{abstract}

\tableofcontents

\section{Introduction}

\subsection{Context} 
Given a variety $X$ over a finite field $k$, its zeta function is the Euler product
\[
\zeta_X(s)=\prod_{x\in X_0} \frac 1 {1-\mathrm{card}(\kappa(x))^{-s}},
\]
where $\kappa(x)$ denotes the residue field at a closed point $x\in X_0$. Since all points have the same residue characteristic $p$, a change of variable $T=p^{-s}$ makes the zeta function $Z_X(T):=\zeta_X(s)$ a power series in the $T$-variable, and it is the exponential of the generating function of point counts of $X$ over the extensions of $\FF_p$. Dwork proved \cite{Dwork1960} that this function was actually rational with coefficients in $\QQ$, and it was shown in the setting of the Weil conjectures that more precisely, the zeta function is the alternating product of the characteristic polynomials of the geometric Frobenius acting on $\ell$-adic cohomology with compact support of the geometric fiber $X_{\overline{k}}$:
\[
Z_X(T)=\prod_{i}\det(1-\phi\cdot T\mid H^i_c(X_{\overline{k}},\QQ_\ell))^{(-1)^{i+1}}.
\]
In particular, the values or more generally the leading coefficient in the Taylor expansion around a given point (called the special value at that point) are given in terms of eigenvalues of the Frobenius acting on the cohomology groups $H^i_c(X_{\overline{k}},\QQ_\ell)$. However, Schneider \cite{Schneider1982} and later Milne \cite{Milne1986} showed that the full knowledge of the eigenvalues is not necessary; namely, the special values at integers $s=n$ are rational numbers so can be understood prime-by-prime, and their $\ell$-part is expressed in terms of $\ell$-adic cohomology of $X$ itself, while for the $p$-part and $n>0$ a further correcting factor, colloquially referred to as Milne's correcting factor, has to be introduced. Those formulas appear slightly unsatisfying since they work prime-by-prime and thus fail to capture the special value on the nose as a rational number. In \cite{Lichtenbaum2005}, Lichtenbaum introduced the Weil-étale cohomology, an integral finitely generated cohomology for smooth projective varieties in characteristic $p$, which is an integral model of $\ell$-adic cohomology. The latter property means that the special value formula at $s=0$ can be written down on the nose in terms of Weil-étale cohomology. Geisser generalized this to integral twists \cite{Geisser2004} and obtained reinterpretations of special value formulas at all integers in terms of Weil-étale cohomology, assuming that the Weil-étale cohomology groups are finitely generated.

On the number theory side, several special value conjectures for $L$-functions of motives appeared, first up to a rational number by Beilinson and Deligne, then integrally by Bloch--Kato \cite{Bloch2007} and Fontaine--Perrin-Riou \cite{Fontaine1991}, and finally integrally for motives with an action of a semi-simple algebra \cite{Burns2001}; the latter is called the equivariant Tamagawa number conjecture (ETNC for short). On the one hand, work of Geisser--Suzuki \cite{Geisser2020} in characteristic $p$ and of the author \cite{Morin2023a} over $\ZZ$, both using Weil-étale cohomology, suggested that a general investigation of special values of $L$-functions associated to $\ZZ$-constructible sheaves on varieties in characteristic $p$ might be fruitful; on the other hand, Burns--Kakde \cite{Burns} proved a special value formula at $s=0$ for non-commutative $L$-functions in characteristic $p$ coming from a Galois cover of varieties over a finite field, suggesting that further investigation of non-commutative $L$-functions in characteristic $p$ was necessary. In the spirit of the ETNC, one should thus consider constructible sheaves of $R$-modules on a variety over a finite field, where $R$ is an order in a semi-simple algebra over $\QQ$, and their associated non-commutative $L$-function, defined analogously to the one in the ETNC. The aim of the present article is to formulate and prove a special values expression at negative integers for those $L$-functions.

\subsection{Results}

\paragraph{Aim of the article}Let $R$ be a $\ZZ$-order in a semi-simple finite-dimensional $\QQ$-algebra $A$ and let $X$ be a scheme separated and of finite type over a finite field. In this article we show, using results of Witte \cite{Witte2014,Witte2016} and the approach of Burns--Kakde \cite{Burns}, that Weil-étale cohomology with compact support, when it is perfect, computes special values at negative integers of non-commutative $L$-functions for $R$-constructible complexes over $X$. We also show that the similar result holds for Weil-eh cohomology with compact support under resolution of singularities in dimension less than $\dim X$.

\paragraph{Non-commutative \texorpdfstring{$L$}{L}-functions} We fix a finite field $k=\FF_q$ of characteristic $p$ and let $f:X\to \Spec(k)$ be a finite type, separated morphism. Let $F$ be a constructible complex of sheaves of $R$-modules on $X$. The semi-simple algebra $A$ has a reduced norm map taking values in the units in its center $\ct{A}$, and similarly for $A[[T]]$, so we can define the non-commutative $L$-function associated to $F$ by the converging Euler product
\begin{defi*}[see \cref{def:L_function}]
	\[
	L_X(F,T):=\prod_{x\in X_0}\mathrm{Nrd}_A(1-\phi_x\cdot T \mid F_{x}\otimes_R A[[T]])^{-1} \in \ct{A}[[T]],
	\]
	where $\phi_x$ denotes the geometric Frobenius at a closed point $x\in X_0$, and $F_{x}$ is the stalk of $F$ at that point.
\end{defi*}
The center of $A$ is a finite product of number fields, so this non-commutative $L$-function can be understood as a family of classical $L$-functions. Writing $\ct{A}\simeq \prod_{i=1}^d K_i$, we can thus talk about values or special values of elements of $\ct{A}[[T]]$ by reasoning component-wise: an element $\eta\in\ct{A}[[T]]$ is said to have a well-defined special value at $x\in \QQ\backslash\set{0}$ if it all its component $\eta_i\in K_i[[T]]$ converge in a small neighborhood of $0$ and extend to meromorphic functions around $x$; the vanishing order $\ord_{t=x}\eta\in \ZZ^d$ and special value $\eta^\ast(x)\in \ct{A}^\times$ is then the collection of vanishing orders and leading coefficient in the Taylor expansion around $x$ (see \cref{def:special_value}). The first observation, which is a consequence of the rationality of the $L$-function proven in our generality by Witte (\cite{Witte2014}, see \cref{thm:rationality}), is that the special value of $L_X(F,T)$ at $t=x\in \QQ\backslash\set{0}$ is always well-defined.

\paragraph{Weil-étale cohomology}  We define the negative Tate twist $F(n)$ of $F$, for $n<0$, by tensoring with the negative tensor power sheaf of all prime-to-$p$ roots of unity. Following Lichtenbaum, we define the Weil-étale cohomology of $X$ with coefficients in $F(n)$ and compact support $R\Gamma_{\mathrm{W,c}}(X,F(n))\in D(R)$  (\cref{def:weil_etale_cohom}). Similarly, following Geisser \cite[def. 5.1]{Geisser2006}, we define the Weil-eh cohomology of $X$ with coefficients in $F(n)$ and compact support $
R\Gamma_{\mathrm{Wh,c}}(X,F(n))$, seeing $F(n)$ as a sheaf for Geisser's $\mathrm{eh}$-topology. Throughout this article, our results will depend on one of the following standard assumptions:
\begin{ass*}[$\mathbf{L}_{F,n}$]
	$R\Gamma_{\mathrm{W,c}}(X,F(n))$ is a bounded complex with finite type cohomology over $\ZZ$.
\end{ass*}
\begin{ass*}[$\mathbf{R}(d)$]
	The strong of form resolution of singularities for varieties up to dimension $d$ over $\overline{\FF}_p$ as in \cite[def. 2.4]{Geisser2006} holds.
\end{ass*}
We investigate some cases where $\mathbf{L}_{F,n}$ can be shown to hold; in particular it holds for locally constant $\ZZ$-constructible sheaves on smooth projective schemes over $k$ (\cref{prop:et_perfect}). Under those assumptions, both complex enjoy nice properties, in particular they are integral models of the $\ell$-adic cohomology of the $\ell$-adic completion of $F(n)$ and their structure rationally forces some semi-simplicity of the geometric Frobenius on the $\ell$-adic cohomology of the base change to $k^{\mathrm{sep}}$.

\paragraph{The Weil-étale Euler characteristic and the special value theorem}

By semi-simplicity of $A$, there is an isomorphism $\mathrm{rrank}_A:K_0(A) \xrightarrow{\simeq} \ZZ^d$ which is called the \emph{reduced rank}. Let $\ct{A}$ denote the center of $A$ and $\cal{O}_{\ct{A}}^\times$ denote the units  in the ring of integers of $\ct{A}$. The natural structures on Weil-étale cohomology, resp. Weil-eh cohomology enable one to define a Weil-étale Euler characteristic $\chi_{W,X}(F(n))\in \faktor{\ct{A}^\times}{\cal{O}_{\ct{A}}^\times}$ (see \cref{def:weil_etale_euler_char}), resp. a Weil-eh Euler characteristic $\chi_{Wh,X}(F(n))$ using Fukaya--Kato's formalism of trivialized determinants (which is recalled in \cref{sec:trivialized_determinant}), and our main theorem is that the special value of $L_X(F(n),T)$ at $t=q^{-n}$, $n\leq 0$, is given up to a unit by the inverse of the Weil-étale Euler characteristic:
\begin{introthm}[see \cref{thm:main} for a more precise statement using Burns--Flach's extended boundary map]\label{thm:introA}
	Let $n\leq 0$.
	\begin{itemize}
		\item Assume $\mathbf{L}_{F,n}$. Then $R\Gamma_{\mathrm{W,c}}(X,F(n))$ is a perfect complex of $R$-modules, we have the vanishing order formula\footnote{If $n<0$, $F(n)$ is torsion so the formula says that the vanishing order is zero}
		\[
		\mathrm{ord}_{t=q^{-n}}L_X(F,T)=\sum_i (-1)^i\cdot i \cdot \mathrm{rrank}_A H^i_{\mathrm{W,c}}(X,F(n))_\QQ
		\]
		and the special value formula
		\[
		L^\ast_X(F,q^{-n})\equiv \frac{1}{\chi_{\mathrm{W},X}(F(n))} \mod \cal{O}_{\ct{A}}^\times.
		\]
		\item Assume $\mathbf{R}(\dim X)$. The similar result holds using $R\Gamma_{\mathrm{Wh,c}}(X,F(n))$ .
	\end{itemize}
\end{introthm}
In particular, the theorem holds unconditionally for locally constant constructible sheaves on smooth projective varieties. This theorem is proven by following Burns--Kakde's strategy of reduction to the $\ell$-adic and $p$-adic case \cite{Burns}. In the local case, the $L$-function is known to be rational by results of Witte \cite{Witte2014} generalizing the Grothendieck--Lefschetz trace formula and is given by the inverse characteristic rational fraction of the geometric Frobenius acting on the cohomology of the $\ell$-adic completion of $F(n)$. This is not quite true at $p$ but again using results of Witte (\cite{Witte2016}, see \cref{thm:witte2}) one can show that the special value is the same under the boundary map $\delta:\ct{A\otimes_\QQ \QQ_\ell}^\times \to K_0(R\otimes_\ZZ \ZZ_\ell,A\otimes_\QQ \QQ_\ell)$ (see \cref{thm:same_vs_under_boundary}).\footnote{Here we crucially need that $n\leq 0$.} Thus the problem reduces to understanding special values of characteristic rational fractions, which is the content of our core theorem, \cref{thm:cohomological_special_value}.
Ultimately, this theorem relies on a nice homological algebra lemma due to Burns (see \cref{lem:burns}). When applied to our situation, the semi-simplicity hypothesis is automatic by the properties of Weil-étale cohomology.

\paragraph{Application to motivic \texorpdfstring{$L$}{L}-functions}
When applying our result to the pushforward of an Artin motive to a smooth projective curve, we obtain a special value formula for non-commutative Artin $L$-functions. Let $K$ be a global field of characteristic $p$, let $S$ be a finite set of places of $K$ and let $V$ be a finite $A$-module with a discrete action of $G_K$, which we identify with an $A$-linear Artin motive over $K$. Denote by $C$ the smooth projective curve with function field $K$ and $U:=C\backslash S$ its open subscheme. Denote by $g:\Spec(K)\to U$ the inclusion of the generic point. 

\begin{introcor}[Special values for non-commutative Artin $L$-functions, see \cref{cor:special_value_artin}]
	There exists a flat $R$-constructible sheaf $F$ on $U$ such that $F\otimes \QQ=g_\ast V$. For such a sheaf, we have for the non-commutative Artin $L$-function $L_S(V,T)$ omitting the places in $S$:
	\[
	L_S(V,T)=L_U(F,T),
	\]
	and thus the vanishing order and special value formulas for any $n\leq 0$:
	\begin{gather*}
		\mathrm{ord}_{t=q^{-n}} L_S(V,T)= \sum_i (-1)^i\cdot i \cdot \mathrm{rrank}_{A} H^i_{\mathrm{W,c}}(U,F(n))_{\QQ},\\
		L_S^\ast(V,q^{-n})\equiv \frac{1}{\chi_{\mathrm{W},U}(F(n))} \mod \cal{O}_{\ct{A}}^\times.
	\end{gather*}
\end{introcor}

The same kind of argument applies to an object $M$ of the derived category of pure $A$-linear motives over $K$ whose non-commutative $L$-function is equal to that of a locally constant $R$-constructible sheaf $F$ on a smooth projective scheme $X$.
\begin{introcor}[Special values for total motives of ``constructible nature'', see \cref{thm:special_value_motive}]
	In the above setting, denote by $F$ the provided locally constant $R$-constructible sheaf on $X$. By definition, we have for the non-commutative motivic $L$-function $L_K(M,T)$ of $M$:
	\[
	L_K(M,T)=L_X(F,T).
	\]
	Hence, we get the vanishing order and special value formulas at $t=q^{-n}$:
	\begin{gather*}
		\mathrm{ord}_{t=q^{-n}} L_K(M,T)= \sum_i (-1)^i\cdot i \cdot \mathrm{rrank}_{A} H^i_{\mathrm{W,c}}(X,F(n))_{\QQ},\\
		L_K^\ast(M,q^{-n})\equiv \frac{1}{\chi_{\mathrm{W},X}(F(n))} \mod \cal{O}_{\ct{A}}^\times.
	\end{gather*}
\end{introcor}
We show that this applies in particular to the case of an unramified Galois extension $L/K$ with group $G$ and a smooth projective variety $X_K$ over $K$ with good reduction, for the induced motive to $K$ of the total derived motive $\oplus h^i(X_L,0)[-i]$ of $X_L$, which is naturally a $\QQ[G]$-linear motive.

\subsection{Overview}
In \cref{sec:background}, we fix definitions and recall the tools that we will use: constructible sheaves, Weil-étale cohomology, relative $K$-theory, Burns--Flach's extended boundary map and Fukaya--Kato's trivialized determinants. We then state precisely the main theorem \cref{thm:main}. In \cref{subsec:artin_induction}, we use Artin induction to prove that locally constant constructible sheaves of $\ZZ$-modules on smooth projective varieties have perfect Weil-étale cohomology. This can be applied to the results of \cref{subsec:weil_etale_perfect}, which say that Weil-étale cohomology is perfect as a complex of $R$-modules as soon as it is as a complex of abelian groups and that Weil-étale cohomology is an integral model of $\ell$-adic cohomology when it is perfect. In \cref{sec:reduction_local_case}, we start the proof of the main theorem, by showing how to reduce its statement to a local statement for characteristic rational functions coming from the geometric Frobenius acting on the cohomology with compact support of an $\ell$-adic sheaf or complex. This works perfectly fine for $\ell\neq p$ by the Grothendieck's cohomological interpretation of $L$-functions of $\ell$-adic sheaves, shown in the necessary generality by Witte. This interpretation breaks down at $p$ but the situation can be salvaged using results of Emerton--Kisin, generalized by Witte. We then investigate in \cref{sec:semisimplicity} a notion of semi-simplicity which enables one to extract informations about special values of characteristic rational functions. In \cref{sec:special_value_cohomological}, the technical core of the article, we prove a special value formula for characteristic rational functions asssociated to an endomorphism of a perfect complex, assuming the adequate semi-simplicity hypothesis. Under this hypothesis, we can reduce the proof to a nice homological algebra lemma \cref{lem:burns} due to Burns, of which we give a new proof. Finally, in \cref{sec:proof} we wrap up the proof of the main theorem and in \cref{sec:applications} we give applications to Artin motives and total derived motives with good reduction over characteristic $p$ global fields.

\subsection{Acknowledgments}

I would like to thank Fabien Pazuki for his support and advice, Matthias Flach, Thomas Geisser, Baptiste Morin and Niranjan Ramachandran for helpful discussions, and the people at the University of Copenhagen, in particular Quingyan Bai, Robert Burklund and Shachar Carmeli, for general discussions around coffee. Special thanks go to Adel Betina for his explanations on $p$-adic meromorphic functions.

\section{Background, notations and statement of the special value theorem}\label{sec:background}

\subsection{Non-commutative \texorpdfstring{$L$}{L}-functions for constructible sheaves}

We fix a finite field $k=\FF_q$ of characteristic $p$ and let $f:X\to \Spec(k)$ be a finite type, separated morphism. Let $R$ be an order in a semi-simple $\QQ$-algebra $A$, that is a finite flat associative $\ZZ$-algebra $R$ such that $A:=R\otimes_\ZZ \QQ$ is a semi-simple (finite-dimensional) $\QQ$-algebra. We will let $\ct{A}$ denote the center of $A$. Let $F$ be an $R$-constructible complex of sheaves on $X$; for us, this will mean a complex of sheaves of right $R$-modules $F\in D(X_{\mathrm{et}},R)$ for which there exists a finite stratification of locally closed subschemes along which $F$ is locally constant with values in $\Perf(R)$, the derived $\infty$-category of perfect complexes of left $R$-modules. Equivalently, $F$ is a bounded complex with constructible cohomology sheaves (in the finite type sense) and finite tor-dimension, see \cite[Rmk. 7.2]{Hemo2023}. This defines a full stable subcategory $\Cons(X,R)$ of $D(X_{\mathrm{et}},R)$.\footnote{Also denoted $D^b_{\mathrm{ctf}}(X,R)$ in \cite[Prop.-def. 4.6]{sga412rapport}} For a closed point $x\in X_0$, the geometric Frobenius\footnote{That is, the element $(\alpha \mapsto \alpha^{p^{-[\kappa(x):\FF_p]}})\in \mathrm{Gal}(\kappa(x)^{\mathrm{sep}}/\kappa(x))$.} $\phi_x$ at $x$ acts on the perfect complex of $R$-modules $F_{\overline{x}}$, which defines an automorphism $1-\phi_x\cdot T$ of $F_{\overline{x}}\otimes_R R[[T]]$ and hence a class $\langle 1-\phi_x\cdot T \mid F_{\overline{x}}\otimes_R R[[T]]\rangle\in K_1(R[[T]])$. The product over all closed points converges for the natural topology\footnote{Another natural topology would be the adic topology given by the Jacobson radical $J(R)$ of $R$; however $J(R)$ is nilpotent (\cref{lem:jabson_order}) and $J(R[[T]])=(T,J(R))$, so the pro-objects $\{R[T]/T^n\}$ and $\{R[[T]]/J(R[[T]])^n\}$ are equivalent.} on $\widehat{K}_1(R[[T]]):=\lim K_1(R[T]/T^n)$ to define the $L$-function element $\mathscr{L}_X(F)$:
\begin{defi}\label{def:L_function}
	\[
	\mathscr{L}_X(F):=\prod_{x\in X_0} \langle 1-\phi_x\cdot T \mid F_{\overline{x}}\otimes_R R[[T]]\rangle^{-1} \in \widehat{K}_1(R[[T]])
	\]
\end{defi}
Since $A$ is semi-simple, we have $A\simeq \prod_{i=1}^d M_{k_i}(D_i)$ where each $D_i$ is a division algebra, and we can write its center $\ct{A}$ as
\[
\ct{A} \xrightarrow{\simeq} \prod_{i=1}^n K_i
\]
where $d=\#(\Spec\ct{A})$ and the $K_i=\ct{D_i}$ are number fields. The reduced norm map induces a map $\widehat{K}_1(A[[T]])\xrightarrow{\mathrm{Nrd}_A} \ct{A}[[T]]^\times$ and we put
\begin{defi}
\[
L_X(F,T):=\mathrm{Nrd}_A(\mathscr{L}_X(F)\otimes_{R[[T]]}A[[T]]) \in \ct{A}[[T]]^\times
\]
\end{defi}

Special values for power series over $\ct{A}$ are defined componentwise:
\begin{defi*}\label{def:special_value}
	Let $\eta\in \ct{A}[[T]]$ and let $x\in \QQ\backslash\set{0}$. We say that $\eta$ has a well-defined special value at $t=x$ if for all $i$, the component $\eta_i$ of $\eta$ in $K_i[[T]]$ extends to a meromorphic function around $x$ (under all embeddings $K_i\hookrightarrow \CC$), with vanishing order $r_i$ (independant of the embedding) and leading term\footnote{We choose the convention that the leading term is of the form $(x-t)^k\cdot \alpha $ in order to avoid introducing signs when studying characteristic polynomials of the form $\det(1-\theta\cdot t \mid V)$} $\eta_i^\ast(x):=\lim_{t\to x} (x-t)^{-r_i}\eta_i(t)\in \CC^\times$ such that $\eta_i^\ast(x)\in K_i^\times$ is independant of the embedding.
	
	If that is the case, we denote by $\eta^\ast(x)=(\eta_i^\ast(x))_i\in \ct{A}^\times$ the corresponding element, called the special value of $\eta$ at $t=x$, and we call the family $\mathrm{ord}_{t=x}\eta:=(r_i)\in \ZZ^d$ the vanishing order of $\eta$ at $x$.
\end{defi*}

\subsection{Weil-étale cohomology}

\paragraph{} Fix for the remainder of the article a compactification $f=gj$ with $j:X\to Y$ an open immersion and $g:Y\to\Spec(k)$ proper of finite type.

\begin{defi}
For $n<0$, we let
\begin{gather*}
\ZZ(n):=\QQ/\ZZ(n)[-1]:=\colim_{p\nmid k} \mu_k^{\otimes n}[-1]\\
F(n):= F\otimes \ZZ(n)
\end{gather*}
and we fix by convention $F(0):=F$. We call $F(n)$ the $n$-th motivic twist of $F$.
\end{defi}

\begin{rmk}
We do not attempt to define positive Tate twist since the naive definition as a tensor product with Bloch's cycle complex does not work well with pushforwards, even finite ones\footnote{Consider the closed immersion $i$ of a closed point $x$ inside a smooth projective curve $C$ over a characteristic $p$ field. Then $(i_\ast \ZZ)\otimes \ZZ(1)=i_\ast \mathcal{O}_{C,x}^{\mathrm{sh}}[-1]$; but for a prime $\ell\neq p$, $i_{\mathrm{proet},\ast} (\ZZ_\ell(1))= (i_{\mathrm{proet},\ast} \ZZ_\ell)\otimes \ZZ_\ell(1)$ since $\ZZ_\ell(1)=\lim \mu_{\ell^n}$, $i^\ast \mu_{\ell^n}=\mu_{\ell^n}$ by étaleness, and $i_{\mathrm{proet}}^\ast$ commutes with all limits.}, and furthermore the methods in this article to analyze the $p$-part of special value at $q^{-n}$ do not work when $n>0$.
\end{rmk}

\begin{defi}\label{def:weil_etale_cohom}
	Let $n\leq 0$, let $W\subseteq \Gal(k^{\mathrm{sep}}/k)$ denote the Weil group of $k$, that is the subgroup generated by the geometric Frobenius $\phi:x\mapsto x^{1/q}$, and let $Y_{\overline{k}}:=Y\times_{\Spec(k)} \Spec(\overline{k})$ denote the geometric fiber of $Y$.
	
	\begin{itemize}
		\item  The pullback of $j_!F(n)$ to $Y_{\overline{k}}$ is a $\Gal(k^{\mathrm{sep}}/k)$-equivariant sheaf, hence by restriction also a $W$-equivariant sheaf, and we define the Weil-étale cohomology of $X$ with coefficients in $F(n)$ as the object
		\begin{align*}
		R\Gamma_{\mathrm{W,c}}(X,F(n))&:=R\Gamma(W,R\Gamma(Y_{\overline{k}},j_!F(n)))\\
		&\simeq \fib\left(R\Gamma(Y_{\overline{k}},j_!F(n))\xrightarrow{1-\phi} R\Gamma(Y_{\overline{k}},j_!F(n))\right)\in D(R)
		\end{align*}
		\item We can consider $F(n)$ as a sheaf for Geisser's $\mathrm{eh}$-topology \cite{Geisser2006}.\footnote{We mean by this the left derived functor of pullback applied to $F$. The pullback functor is a priori not left $t$-exact, but it is if one assumes resolution of singularities.} We define the Weil-eh cohomology of $X$ with coefficients in $F(n)$, following \cite[def. 5.1]{Geisser2006}, as
		\[
		R\Gamma_{\mathrm{Wh,c}}(X,F(n)):=R\Gamma(W,R\Gamma_{\mathrm{c}}((X_{\overline{k}})_{\mathrm{eh}},F(n)))\in D(R)
		\]
	\end{itemize}
\end{defi}

\begin{rmk}
The first notation is evil: it doesn't reflect the fact that it very much depends on the chosen compactification, but it reduces symbols overload. The second definition corrects this defect, but to say anything about it one has to usually assume resolution of singularities over $\overline{\FF}_p$ for schemes of dimension less than $\dim X$.	
\end{rmk}

\paragraph{} We will rely on the following standard assumptions to get our results:
\begin{ass*}[$\mathbf{L}_{F,n}$]
	$R\Gamma_{\mathrm{W,c}}(X,F(n))$ is a bounded complex with finite type cohomology over $\ZZ$.
\end{ass*}
\begin{ass*}[$\mathbf{R}(d)$]
	The strong of form resolution of singularities for varieties up to dimension $d$ over $\overline{\FF}_p$ as in \cite[def. 2.4]{Geisser2006} holds.
\end{ass*}

By Artin-Schreier theory, $H^2_W(\mathbb{A}^1_k,\ZZ)$ is not finitely generated \cite[p. 685]{Geisser2004}, hence compactifying is necessary. Moreover, using the étale site is also problematic in the proper non-smooth case: Geisser gives an example \cite[Ex. 3]{Geisser2006} of a normal projective surface $X$ with one singular point $P$ and smooth projective blowup $Y$ at $P$ with exceptional divisor a nodal projective curve, such that $H^3_W(X,\ZZ)$ is not finitely generated. Moreover, denoting by $U=X\backslash P\simeq Y\backslash C$ and  $j:U \to X$ and $j':U \to Y$ the open immersions, we have $H^2_W(X,j_!\ZZ)\otimes \QQ=0$ while $H^2_W(Y,j'_!\ZZ)\otimes \QQ=\QQ$, so Weil-étale cohomology with compact support of $U$ depends on the choice of compactification.

The fix to the above issues is to use Weil-eh cohomology, but it necessitates the use of resolution of singularities (see \cref{prop:eh_perfect}). We mention that $R(2)$ is known by Abhyankar, as well as $R(3)$ for $p\geq 7$. By work of Cossart and Piltant, the part of $R(3)$ about resolution of singularities is known in general, however I do not know whether the part about refining a proper birational morphism from a smooth scheme as a succession of blowups with smooth centers has been shown.
%

\subsection{Trivialized determinants}\label{sec:trivialized_determinant}

We recall Fukaya--Kato's formalism of trivialized determinants: given a ring map $R\to S$, consider the Picard groupoids $V_R$ and $V_S$ attached to the universal determinants $[-]_R:\Perf(R)^\simeq\to V_R$, $[-]_S: \Perf(S)^\simeq\to V_S$ for $R$ and $S$; their $\pi_i$, $i=0,1$, compute the lower $K$-groups of $R$ and $S$ respectively \cite[§2.3]{Burns2001}, \cite{Deligne1987}. A class in $K_0(R)$, written as $[P]-[Q]$ where $P$ and $Q$ are projective modules of finite type, corresponds to the isomorphism class of $[P]_R[Q]_R^{-1}$ in $\pi_0(V_R)$, while  if $\psi$ is an automorphism of $C\in \Perf(S)$, the class $\langle \psi \mid C\rangle\in K_1(S)$ is represented by the composite 
$\begin{tikzcd}
	{\mathbf{1}_S} & {[C]_S[C]_S^{-1}} & {[C]_S[C]_S^{-1}} & {\mathbf{1}_S}
	\arrow["{\mathrm{can}^{-1}}", "\simeq"',from=1-1, to=1-2]
	\arrow["{[\psi]\cdot\mathrm{id}}", "\simeq"', from=1-2, to=1-3]
	\arrow["{\mathrm{can}}",  "\simeq"',from=1-3, to=1-4]
\end{tikzcd}$ in $\pi_1(V_S):=\mathrm{Aut}(\mathbf{1}_S)$ (see \cite[§2.3]{Burns2001}).

Define $V_0$ as the Picard groupoid with one object and the identity automorphism; then we can form the pullback in categories
\[\begin{tikzcd}[ampersand replacement=\&]
	{V(R,S)} \& {V_0} \\
	{V_R} \& {V_S}
	\arrow[from=1-1, to=1-2]
	\arrow[from=1-1, to=2-1]
	\arrow["{\usebox\pullback}"{description, pos=0}, draw=none, from=1-1, to=2-2]
	\arrow[from=1-2, to=2-2]
	\arrow[from=2-1, to=2-2]
\end{tikzcd}\]
which is a Picard groupoid. One can show that $\pi_0 V(R,S)$ computes the relative $K$-group $K_0(R,S)$ \cite[adapt prop. 2.5]{Burns2001}. Since the objects of $V(R,S)$ can be described explicitly as pairs $(X,\alpha)$ where $X\in V_R$ and $\alpha:X\otimes_R S \xrightarrow{\simeq} \mathbf{1}_S$ is an isomorphism, such a pair defines functorially a class $\chi(X,\alpha)$ in $K_0(R,S)$ called the refined Euler characteristic. We also observe that the map $K_1(S)\to K_0(R,S)$ is modeled by $\alpha\in \mathrm{Aut}(\mathbf{1}_S) \mapsto (\mathbf{1}_R,\alpha)$ and the map $K_0(R,S)\to K_0(R)$ by the projection $(X,\alpha)\mapsto X$.

\subsection{The extended boundary map}
In this short section we recall the definition of the extended boundary map $\widehat{\delta_R}:\ct{A}^\times \to K_0(R,A)$, and record some observations on its kernel. The special value theorem will determine special values of the non-commutative $L$-function up to an element of that kernel. The main takeaway is that elements in the image of $\widehat{\delta_R}$ are determined, at worst, up to a unit in the ring of integers of $\ct{A}$.

We denote by $R_\ell$ the base change $R\otimes_\ZZ \ZZ_\ell=R^\wedge_\ell$ of $R$ to $\ZZ_\ell$ and $A_{\ell}:=R_\ell[\frac 1 \ell]=A\otimes_\QQ \QQ_{\ell}$. The reduced norm for $A_{\ell}$ is an isomorphism \cite[Thm. 45.3]{Curtis1987}, and we define the boundary map $\delta_{\ell}$ as the composite
$\delta_\ell: \ct{A_\ell}^\times \xrightarrow[\simeq]{\mathrm{Nrd}_{A_{\ell}}^{-1}} K_1(A_{\ell}) \to K_0(R_\ell,A_{\ell})$.
There is a factorization \cite[Lem. 9]{Burns2001}
\[\begin{tikzcd}
	{\prod_\ell\ct{A_\ell}^\times} & {\prod_\ell K_0(R_{\ell},A_{\ell})} & {} \\
	& {\bigoplus_\ell K_0(R_{\ell},A_{\ell})} \\
	{\ct{A}^\times} & {K_0(R,A)}
	\arrow["{(\delta_{\ell})}", from=1-1, to=1-2]
	\arrow["{\widehat{\delta_R}}"', dashed, from=3-1, to=3-2]
	\arrow[from=3-1, to=1-1]
	\arrow["\simeq", from=3-2, to=2-2]
	\arrow[hook, from=2-2, to=1-2]
\end{tikzcd}\]
giving the \emph{extended boundary map} $\widehat{\delta_R}:\ct{A}^\times \to K_0(R,A)$. The $K$-group $K_1(A)$ is a subgroup of $\ct{A}^\times$ via the reduced norm \cite[Thm. 45.3]{Curtis1987} and the restriction $\widehat{\delta_R}\circ \mathrm{Nrd}_{A}$ to $K_1(A)$ is the canonical boundary map $K_1(A)\to  K_0(R,A)$.

\begin{prop}\label{prop:kernel_extended_boundary}
	\begin{enumerate}
		\item[]
		\item We have 
		\[\ker(\widehat{\delta_R})=\set{\alpha\in \ct{A}^\times \mid \forall \ell,~ \alpha \in \mathrm{Nrd}_{A_\ell}(K_1(R_\ell))\subseteq\ct{A_\ell}^\times}.
		\]
		\item When $R$ is commutative, we have $\ker(\widehat{\delta_R})=R^\times$.
		\item When $R$ is an hereditary order in a semi-simple algebra with center $\ct{A}=\prod_{i=1}^d K_i$ where the $K_i$ are number fields, we have $\ker(\widehat{\delta_R})=\mathcal{O}_{\ct{A}}^\times:=\prod_{i=1}^d\cal{O}_{K_i}^\times$.
		\item In general, we have $\ker(\widehat{\delta_R})\subseteq \mathcal{O}_{\ct{A}}^\times$, and the quotient $\ct{A}^\times/\mathcal{O}_{\ct{A}}^\times$ is torsion-free.
	\end{enumerate}
\end{prop}

\begin{proof}
	\begin{enumerate}
		\item This comes from the identifications $K_1(A_\ell)\simeq \ct{A_\ell}^\times$, induced by the reduced norm, and $K_0(R,A)\simeq \bigoplus_\ell K_0(R_\ell,A_\ell)\subseteq \prod_\ell K_0(R_\ell,A_\ell)$.
		\item Since $R_\ell$ is a semi-local commutative ring, we have $K_1(R_\ell)\simeq R_\ell^\times$, so this is a consequence of the pullback diagram of abelian groups and of rings, usually called the arithmetic fracture square
		\[\begin{tikzcd}[ampersand replacement=\&]
			R \& {\prod_\ell R_\ell} \\
			A \& {(\prod_\ell R_\ell)\otimes_\ZZ \QQ.}
			\arrow[from=1-1, to=1-2]
			\arrow[from=1-1, to=2-1]
			\arrow["{\usebox\pullback}"{description, pos=0}, draw=none, from=1-1, to=2-2]
			\arrow[from=1-2, to=2-2]
			\arrow[from=2-1, to=2-2]
		\end{tikzcd}\]
		\item If $R$ is hereditary, the $\ell$-adic completions $R_\ell$ are again hereditary. If $\ct{A_\ell}\simeq \prod_{j=1}^f L_j$ is a decomposition into local fields, put $\mathcal{O}_{\ct{A_\ell}}:= \prod_{j=1}^f \mathcal{O}_{L_j}$, and define similarly $\mathcal{O}_{\ct{A}\otimes_\QQ \QQ_\ell}$. By \cite[Thm. 2]{Wilson1977} we have $\mathrm{Nrd}_{A_\ell}(K_1(R_\ell))=\mathcal{O}_{\ct{A_\ell}}^\times$. By reasoning with valuations, one observes that  $(\ct{A}\otimes_\QQ \QQ_\ell)\cap \mathcal{O}_{\ct{A_\ell}}=\mathcal{O}_{\ct{A}\otimes_\QQ \QQ_\ell}=(\mathcal{O}_{\ct{A}})_\ell$, whence the result again by the arithmetic fracture square.
		\item Any order embeds into a maximal one, which must be hereditary \cite[Thm. 26.12]{Curtis1987}. Moreover, $\ct{A}$ decomposes as a product of number fields, and for any number field $K$, there is an injection
		\[
		K^\times/\mathcal{O}_{K}^\times \hookrightarrow \prod_{\ell} K_\ell^\times/\mathcal{O}_{K_\ell}^\times
		\]
		and the right hand side is a product of free abelian groups of finite rank.
	\end{enumerate}
\end{proof}

\subsection{The Weil-étale Euler characteristic and the special value theorem}

There is a canonical class $e\in H^1_W(\Spec(k),\ZZ)=\Hom(W,\ZZ)$ given by $(\varphi:x\mapsto x^q)\mapsto 1$. This class is represented by the non-split extension of $\ZZ[W]$-modules
\[
0\to \ZZ\to \ZZ^2\to \ZZ\to 0,
\]
where $\varphi$ acts as $\begin{pmatrix}	1 & 1\\ 0 & 1 \end{pmatrix}$. The extension determines a map $\beta:\ZZ \to \ZZ[1]$ whose homotopy class is $e$, and tensoring with the pullback of $\beta$ to $X$ gives a map $F\to F[1]$ of W-equivariant sheaves inducing a map $R\Gamma_{\mathrm{W,c}}(X,F)\to R\Gamma_{\mathrm{W,c}}(X,F)[1]$, which on the level of cohomology groups is cup product with the pullback of $e$. Since $H^2_W(\Spec(k),\ZZ)=0$, we have $(\cup e)^2=0$ so that $(H^\ast_{\mathrm{W,c}}(X,G),\cup e)$ a cochain complex. Rationally, the map $\beta$ is well-understood:

\begin{prop}[\cite{Geisser2004,Geisser2006}]
	Let $G$ be an étale sheaf on $X$ with torsion-free cohomology sheaves. There is a canonical direct sum decomposition
	\[
	R\Gamma_{\mathrm{W,c}}(X,G)\simeq R\Gamma_{\mathrm{et}}(Y,j_! G)[-1] \oplus R\Gamma_{\mathrm{et}}(Y,j_! G)
	\]
	under which the map $\beta$ is identified with the matrix $\begin{pmatrix}
		0 & 1\\ 0 & 0
	\end{pmatrix}$.
	In particular, the complex $(H^\ast_{\mathrm{W,c}}(X,G),\cup e)$ is acyclic. The same holds for Weil-eh cohomology.
\end{prop}

 Using the above formalism of trivialized determinants, one defines a canonical element $\chi_{\mathrm{W},X}(F(n))\in K_0(R,A)$ when $R\Gamma_{\mathrm{W,c}}(X,F(n))$ is perfect:
\begin{defi}\label{def:weil_etale_euler_char}
	\begin{itemize}
		\item[]		
		\item Assume that $R\Gamma_{\mathrm{W,c}}(X,F(n))\in \Perf(R)$. Let $\alpha^\mathrm{W}$ be the isomorphism of $A$-determinants
		\begin{align*}\alpha^\mathrm{W}:[R\Gamma_{\mathrm{W,c}}(X,F(n))_\QQ]_A & \xrightarrow{\simeq} [R\Gamma_{\mathrm{et}}(Y,j_!F(n))_\QQ[-1] \oplus R\Gamma_{\mathrm{et}}(Y,j_!F(n))_\QQ
			]_A\\
			&\xrightarrow{\simeq} [R\Gamma_{\mathrm{et}}(Y,j_!F(n))_\QQ]_A^{-1}[R\Gamma_{\mathrm{et}}(Y,j_!F(n))_\QQ]_A\\
			& \xrightarrow[\simeq]{\mathrm{can}} \mathbf{1}_A.
		\end{align*}
		We define the Weil-étale Euler characteristic of $F(n)$ as (\emph{cf} \cref{sec:trivialized_determinant})
		\[
		\chi_{\mathrm{W},X}(F(n)):=\chi\left([R\Gamma_{\mathrm{W,c}}(X,F(n))]_R,\alpha^\mathrm{W}\right)\in K_0(R,A).
		\]
		\item Assume that $R\Gamma_{\mathrm{Wh,c}}(X,F(n))\in \Perf(R)$.
		Let $\alpha^{\mathrm{Wh}}$ be the isomorphism of $A$-determinants
		\begin{align*}\alpha^{\mathrm{Wh}}:[R\Gamma_{\mathrm{Wh,c}}(X,F(n))_\QQ]_A & \xrightarrow{\simeq} [R\Gamma_{\mathrm{eh},c}(X,F(n))_\QQ[-1] \oplus R\Gamma_{\mathrm{eh}}(X,F(n))_\QQ
			]_A\\
			&\xrightarrow{\simeq} [R\Gamma_{\mathrm{eh},c}(X,F(n))_\QQ]_A^{-1}[R\Gamma_{\mathrm{eh},c}(X,F(n))_\QQ]_A\\
			& \xrightarrow[\simeq]{\mathrm{can}} \mathbf{1}_A.
		\end{align*}
		We define the Weil-eh Euler characteristic of $F(n)$ as
		\[
		\chi_{\mathrm{Wh},X}(F(n)):=\chi\left([R\Gamma_{\mathrm{Wh,c}}(X,F(n))]_R,\alpha^{\mathrm{Wh}}\right)\in K_0(R,A).
		\]
	\end{itemize}
\end{defi}

By Morita invariance and commutation of $K_0$ with finite products of rings, we find that the canonical composite map
\[
\mathrm{rrank}_A:K_0(A)\xrightarrow{\simeq}\prod_{i=1}^d K_0(D_i) \xrightarrow{(\mathrm{rank}_{D_i})_i} \ZZ^d
\]
which we call the \emph{reduced rank}, is an isomorphism.

\begin{thm}\label{thm:main}
	Let $n\leq 0$.
	\begin{itemize}
		\item Assume $\mathbf{L}_{F,n}$. Then $R\Gamma_{\mathrm{W,c}}(X,F(n))\in\Perf(R)$, we have the vanishing order formula\footnote{If $n<0$, $F(n)$ is torsion so the formula says that the vanishing order is zero}
		\[
		\mathrm{ord}_{t=q^{-n}}L_X(F,T)=\sum_i (-1)^i\cdot i \cdot \mathrm{rrank}_A H^i_{\mathrm{W,c}}(X,F(n))_\QQ
		\]
		and the special value formula
		\[
		\widehat{\delta_R} (L^\ast_X(F,q^{-n}))=-\chi_{\mathrm{W},X}(F(n)).
		\]
		\item Assume $\mathbf{R}(\dim X)$. The similar result holds using $R\Gamma_{\mathrm{Wh,c}}(X,F(n))$ .
	\end{itemize}
\end{thm}
The full proof will be given in \cref{sec:proof}.

\begin{rmk}
	\begin{enumerate}
		\item[]
		\item When $R=\ZZ$, we have $K_1(\QQ)=\QQ^\times$, $K_0(\ZZ,\QQ)=\QQ^\times/\ZZ^\times$ and the boundary map is the reduction map. In that case the above theorem determines the special value up to sign. In general, $K_0(R,A)$ will be more complicated as $K_1(A)\subsetneq \ct{A}^\times$ and $K_0(R,A)\to K_0(R)$ can be non-zero. The class  $\widehat{\delta_R}(L_X^\ast(F,q^{-n}))\in K_0(R,A)$ does not necessarily belong to the image of the boundary map $K_1(A)\to K_0(R,A)$. If it does, by the special value theorem we obtain that $[R\Gamma_{W,c}(X,F(n))]=0$ in $K_0(R)$, which by the argument of \cite[end of §6.2.2]{Burns} implies that the complex $R\Gamma_{W,c}(X,F(n))$ is represented by a finite complex with \emph{free} terms of finite rank. However, $\widehat{\delta_R}(L_X^\ast(F,q^{-n}))$ always belongs to the locally free class group $\mathrm{Cl}(R,A):=\ker(K_0(R,A)\to \prod_\ell K_0(A_\ell))$, by construction.
		\item As we observed in \cref{prop:kernel_extended_boundary}, since the map $R\to A$ must factor through a maximal order and thus $\ct{A}^\times/\ker(\widehat{\delta_R})$ always admits the torsion-free quotient $\ct{A}^\times/\mathcal{O}_{\ct{A}}^\times$; the special value theorem determines the special value, at worst, in this quotient.
		\item It might be possible to obtain special value formulas at $t=x\in \ZZ\backslash\set{0}$ if one manages to define the correct Weil-étale cohomology $R\Gamma_{\mathrm{W,c}}(X,F\set{x})\in \Perf(R)$ such that
		\[
		R\Gamma_{\mathrm{W,c}}(X,F\set{x})\otimes_\ZZ \ZZ_\ell \simeq \fib\left(R\Gamma_{\mathrm{c}}(X_{\overline{k}},F_\ell)\xrightarrow{1-\phi\cdot x}R\Gamma_{\mathrm{c}}(X_{\overline{k}},F_\ell\right)\in \Perf(R_\ell)
		\]
		for all primes $\ell$, and such that its rational structure ``forces $1-\phi\cdot x$ to be semi-simple at $0$'' (see \cref{cor:semisimplicity_frobenius}). One way one might be able to do this is to study a relaxed Weil-étale topos allowing sheaves equivariant for the monoid generated by the geometric Frobenius.
	\end{enumerate}
\end{rmk}

\section{Properties of Weil-étale cohomology}

\subsection{Artin induction and dévissage}\label{subsec:artin_induction}
In this section, we use Artin induction to deduce dévissage results for constructible sheaves. Artin induction enables one to reduce from locally constant sheaves to constant sheaves. We use our dévissage result to show that Weil-étale cohomology, resp. Weil-eh cohomology is perfect over $\ZZ$ for locally constant $\ZZ$-constructible sheaves sheaves on smooth projective schemes (\cref{prop:et_perfect}), resp. for $\ZZ$-constructible sheaves on all finite type separated schemes over $k$ under resolution of singularities (\cref{prop:eh_perfect}).

Let $d\in \NN$. We will denote by $\mathrm{Sch}_k^d$ the category of finite type separated schemes over $k$ of dimension less than $d$. The following dévissage lemma, involving Artin induction and classical arguments, reduces claims about constructible sheaves to claims about constant sheaves on smooth projective varieties.
\begin{lem}[Dévissage lemma]\label{lem:devissage}
	Let $P$ be a property of pairs $(X,F)$ where $X \in \mathrm{Sch}_k^d$ and $F\in \Cons(X,\ZZ)$.
	\begin{enumerate} 
		\item Assume that
		\begin{enumerate}
			\item $P$ holds for any torsion constructible sheaf $F$ (concentrated in degree $0$) and for its shift $F[-1]$;
			\item $P$ holds for pairs $(X,\ZZ)$ where $X$ is smooth projective (resp. smooth and proper) and $\ZZ$ is the constant sheaf;
			\item In a fiber sequence $F\to G\to H$ of $\ZZ$-constructible complexes on $X\in \mathrm{Sch}_k^d$, if $P(X,-)$ holds for $F$ and $H$ then it holds for $G$;
			\item for $f:X\to Y$ is a finite étale morphism, then $P(X,F)$ holds if and only if $P(Y,f_\ast F)$ holds;
			\item $P$ is stable under retracts of complexes.
		\end{enumerate} 
		Then $P$ holds for all pairs $(X,F)$ where $X$ is smooth projective (resp. smooth proper) of dimension less than $d$ and $F$ is a $\ZZ$-constructible \emph{sheaf} over $X$ with locally constant torsion-free part. If moreover $P$ is stable under shifts of complexes then $P$ holds for all $F\in \Cons(X,\ZZ)$ such that its cohomology sheaves have locally constant torsion-free part.
		\item Assume that
		\begin{enumerate}
			\item $P$ holds for any torsion constructible sheaf $F$ (concentrated in degree $0$) and for its shift $F[-1]$;
			\item $P$ holds for pairs $(X,\ZZ)$ where $X$ is smooth projective and $\ZZ$ is the constant sheaf;
			\item In a fiber sequence $F\to G\to H$ of $\ZZ$-constructible complexes on $X\in \mathrm{Sch}_k^d$, if $P(X,-)$ holds for two of the three then it holds for the third;
			\item for $f:X\to Y$ is a quasi-finite morphism, then $P(X,F)$ holds if and only if $P(Y,f_!F)$ holds\footnote{The functor $f_!$ is well-defined in this case.};
			\item $P$ is stable under retracts of complexes;
			\item $\mathbf{R}(d)$ holds.
		\end{enumerate}
		Then $P$ holds for all pairs $(X,F)$ where $X \in \mathrm{Sch}_k^d$ and $F$ is a $\ZZ$-constructible \emph{sheaf} over $X$. If moreover $P$ is stable under shifts of complexes then $P$ holds for all pairs $(X,F)$ where $X \in \mathrm{Sch}_k^d$ and $F\in \Cons(X,\ZZ)$.
	\end{enumerate}
\end{lem}

\begin{proof}
	\begin{enumerate}
		\item[]
		\item As $F$ is the direct sum of its torsion and torsion-free parts, we can assume by $(a)$ and $(c)$ that $F$ is torsion-free. Then $F$ corresponds to a representation of $\pi_1^{\mathrm{proet}}(X)$ in free finite type $\ZZ$-modules. As $X$ is normal, $\pi_1^{\mathrm{proet}}(X)=\pi_1^{\mathrm{et}}(X)$ is profinite and the action factors through a finite quotient. By Artin induction \cite[Prop. 4.1, Cor. 4.4]{Swan1960}, there exists finitely many finite étale morphisms $\pi_i:Y_i\to X,\pi_j:Y_j\to X$ and a fiber sequence
		\[
		N[-1] \to F^n\oplus\bigoplus_i \pi_{i,\ast} \ZZ \to \bigoplus_j \pi_{j,\ast} \ZZ
		\]
		where $N$ is a torsion constructible sheaf. Hence $P(X,F)$ holds by applying all hypotheses. If $P$ is stable under shifts, we can use the standard truncation filtration and $(c)$ to conclude. 
		\item If we put $Q(X):=P(X,\ZZ)$, then $(c)$ and $(d)$ applied to the fiber sequence
		\[
		j_!\ZZ \to \ZZ \to i_\ast \ZZ
		\]
		for an open-closed decomposition $U\xrightarrow{j} X \xleftarrow{i} Z$ shows that $Q$ satisfies two out of three for open-closed decomposition. Thus, using $(f)$ and $(b)$ with \cite[Lem. 2.7]{Geisser2006}, we get that $P$ holds for pairs $(X,\ZZ)$ where $X\in \mathrm{Sch}_k^d$. By induction on the dimension of $X$, let us prove that $P$ holds for any $\ZZ$-constructible sheaf $F$. We can assume that $X$ is reduced. Consider the normalization $\pi:Y\to X$ of $X$; it is finite birational, so the cokernel $Q$ in the short exact sequence
		\[
		0 \to F \to \pi_\ast \pi^\ast F \to Q \to 0
		\]
		is supported on a closed subscheme of dimension $<\dim X$. By the induction hypothesis, $(c)$ and $(d)$ we deduce that if $P(Y,\pi^\ast F)$ holds then $P(X,F)$ holds. Moreover, the finitely many irreducible components of $Y$ are disjoint by construction, so $P(Y,\pi^\ast F)$ holds if and only if $P(-,\pi^\ast F)$ holds on each irreducible component. We can thus assume that $X$ is normal and irreducible. Since $F$ is locally constant on a dense open, we can assume by $(c)$ and the induction hypothesis that $F$ is locally constant. We can now conclude as in $1.$
		
		If $P$ is stable under shifts of complexes, then we can represent $F\in \Cons(X,\ZZ)$ by a bounded complex of torsion-free $\ZZ$-constructible sheaves and use the naive truncation filtration together with $(d)$ to conclude.\qedhere
	\end{enumerate}
\end{proof}

\begin{prop}\label{prop:et_perfect}
	$\mathbf{L}_{F,n}$ holds for $n\leq 0$ and
	\begin{enumerate}
		\item $F$ a locally constant $\ZZ$-constructible complex on a smooth projective scheme (or $n=0$ and a smooth proper scheme).
		\item $F$ a $\ZZ$-constructible complex on a curve;
		\item $F$ a generically unramified\footnote{That is, such that the restriction to a dense open equals the restriction of a locally constant $\ZZ$-constructible sheaf defined on the whole space.} $\ZZ$-constructible complex on a smooth projective surface and $n=0$;
	\end{enumerate}
\end{prop}
\begin{proof}
	The classical finiteness result for cohomology with compact support of a torsion constructible complex, together with the étale--Weil-étale comparison \cite[Cor. 5.2]{Geisser2004} shows that $\mathbf{L}_{F,n}$ is known for such complexes. By \cref{lem:devissage}.1., to show the first claim it suffices to show it for $F=\ZZ$ on a smooth projective scheme $X$, whence it is \cite[Prop. 9.2]{Geisser2004} (or \cite[Thm 7.3, Prop. 7.4]{Geisser2004} for $n=0$ in the smooth proper case).
	
	For the case of a curve, one can apply the strategy of \cref{lem:devissage}.2. since the normalization will be smooth, and the fact that $\mathbf{L}_{\ZZ,n}$ is known for $n\leq 0$ on a curve \cite[Lem. 3.3]{Lichtenbaum2005},\cite{Geisser2004}[Prop. 9.4]. We obtain the case of a smooth projective surface through a similar strategy and \cite[Thm. 3.3]{Lichtenbaum2005}.
\end{proof}

\begin{prop}\label{prop:eh_perfect}
	Under $\mathbf{R}(d)$, Weil-eh cohomology with compact support is perfect for all $\ZZ$-constructible complexes and their negative twists on separated schemes of finite type of dimension $\leq d$.
\end{prop}
\begin{proof}
	This follows from \cref{lem:devissage}.2. by the classical finiteness of étale cohomology with compact support for torsion constructible complexes, \cite[Thm 7.3, Prop. 7.4, Prop. 9.2]{Geisser2004} and \cite[Thm. 3.6, Thm. 5.2, Cor. 5.5]{Geisser2006}.
\end{proof}

\begin{prop}\label{prop:comparison_etale_eh}
	Weil-eh, Weil-étale, eh and étale cohomology with compact support all agree on torsion constructible complexes and their negative twists on any separated finite type scheme over $k$.
	Under $\mathbf{R}(d)$, Weil-eh cohomology and Weil-étale cohomology agree for all locally constant $\ZZ$-constructible complexes and their negative twists on smooth projective schemes of dimension $\leq d$.
\end{prop}
\begin{proof}
	The first claim is \cite[Cor. 5.2]{Geisser2004} and \cite[Thm. 3.6, Thm. 5.2]{Geisser2006}. For the second claim, there is a canonical comparison map, to which one can apply \cref{lem:devissage}.1., \cite[Cor. 5.5]{Geisser2006} and the first claim.
\end{proof}

\subsection{Perfectness of Weil-étale cohomology}\label{subsec:weil_etale_perfect}

In this section, we observe that if the Weil-étale cohomology of an $R$-constructible complex is perfect over $\ZZ$, then it is perfect over $R$ (\cref{cor:perfect_Z_implies_perfect_R}), and it is an integral model for the $\ell$-adic cohomology of the $\ell$-adic completion of the sheaf (\cref{prop:Weil_integral_model}).

\paragraph{}
\begin{lem}\label{lem:perfect_local}
	Let $M\in D(R)$ such that $M\in \Perf(\ZZ)$ and $M\otimes^L_\ZZ \ZZ/\ell\in\Perf(R/\ell)$ for all primes $\ell$. Then $M\in \Perf(R)$.
\end{lem}
\begin{proof}
	Since $R$ is finite over $\ZZ$ it is left Noetherian. The hypothesis implies that all cohomology groups are finite over $\ZZ$, hence also over $R$, so $R\Gamma_{\mathrm{W,c}}(X,F)$ is almost perfect, that is, bounded and pseudocoherent. Thus, $M_\QQ$ is almost perfect and bounded over the semi-simple algebra $A$ so it is perfect. Let $N$ denote a perfect complex over $R$ such that $N_\QQ\simeq M_\QQ$. Then $N$ is compact in $D(R)$ so there is a map $N\to M$ inducing the isomorphism $N_\QQ\simeq M_\QQ$ after tensoring with $\QQ$. The fiber $F=\fib(N\to M)$ is perfect over $\ZZ$ and torsion, so it is killed by an integer $N>0$. We deduce that $M\otimes \ZZ[1/N]$ is perfect. To prove that $M$ is perfect, it thus remains to check that $M$ has bounded Tor-amplitude. The fact that $M\otimes \ZZ[1/N]$ is perfect gives a uniform bound at $\ell\nmid N$; thus if we can prove that $M\otimes \ZZ_\ell$ is perfect for all primes $\ell$, this will give a bound at the finitely many $\ell\mid N$, concluding the proof. The ring $R_{\ell}$ is $\ell$-adically complete, and $M\otimes\ZZ_\ell=M^\wedge_{\ell}$ is also $\ell$-adically complete by the $\ZZ$-perfectness, so by \cite[\href{https://stacks.math.columbia.edu/tag/09AW}{Tag 09AW}]{stacks-project} we get from $M\otimes\ZZ/\ell \in \Perf(R/\ell)$ that $M\otimes\ZZ_\ell\in \Perf(R_\ell)$.
\end{proof}

\begin{cor}\label{cor:perfect_Z_implies_perfect_R}
	\begin{itemize}
		\item[]
		\item Assume $\mathbf{L}_{F,n}$. Then $R\Gamma_{\mathrm{W,c}}(X,F(n))\in \Perf(R)$.
		\item Assume $\mathbf{R}(\dim X)$. Then $R\Gamma_{\mathrm{Wh,c}}(X,F(n))\in \Perf(R)$.
	\end{itemize}
\end{cor}

\begin{proof}
	We prove the first case, the proof is exactly the same in the second case. By \cref{lem:perfect_local}, it suffices to show that $R\Gamma_{\mathrm{W,c}}(X,F(n)/\ell)\in \Perf(R/\ell)$ for any prime $\ell$. Now for torsion coefficients, Weil-étale cohomology and étale cohomology with compact support coincide (\cref{prop:comparison_etale_eh}), hence $R\Gamma_{\mathrm{W,c}}(X,F(n)/\ell)=R\Gamma_{\mathrm{c}}(X,F(n)/\ell)$ belongs indeed to $\Perf(R/\ell)$: $F(n)/\ell$ is a constructible sheaf over the finite ring $R/\ell$ and we conclude with the stability of $D^b_{\mathrm{ctf}}(-,R/\ell)$ under $(-)_!$ \cite[Thm. 4.9]{sga412rapport}.
\end{proof}

In particular, assuming $\mathbf{L}_{F,n}$, resp. $\mathbf{R}(\dim X)$, the Euler characteristic $\chi_{\mathrm{W},X}(F(n))$, resp. $\chi_{\mathrm{Wh},X}(F(n))$, can always be defined.

\paragraph{} Define $F_\ell(n)$ to be the $\ell$-adic completion (in the proétale topos) of $F(n)$
. Then $F_p(n)=0$ for $n<0$ and we have
\[
\begin{array}{lll}
	\colim_k (\mu_{\ell^k}^{\otimes n} \otimes_\ZZ \ZZ/\ell^m) & \simeq \colim_{k\geq m} \mu_{\ell^m}^{\otimes n} \otimes_\ZZ  \ZZ/\ell^k	&\simeq \mu_{\ell^m}^{\otimes n} \otimes_{\ZZ} (\colim_{k\geq m} \ZZ/\ell^k)\\
	&&\simeq   \mu_{\ell^m}^{\otimes n}[1]
\end{array}
\]
hence
\[
\begin{array}{ll}
	F_\ell(n) \simeq \lim_m (F\otimes_\ZZ \colim_k \mu_{\ell^k}^{\otimes n} [-1])/\ell^m &\simeq \lim_m F/\ell^m \otimes_{\ZZ/\ell^m} \colim_k (\mu_{\ell^k}^{\otimes n} [-1]\otimes_\ZZ \ZZ/\ell^m)\\
	&\simeq \lim_m F/\ell^m \otimes_{\ZZ/\ell^m} \mu_{\ell^m}^{\otimes n}\\
	&\simeq \lim_m F_\ell/\ell^m \otimes_{\ZZ_\ell} \ZZ_\ell(n)\\
	& \simeq F_\ell\otimes_{\ZZ_\ell}\ZZ_\ell(n)
\end{array}
\]
where the last equality holds since $\ZZ_\ell(n)$ is dualizable, which implies that $-\otimes_{\ZZ_\ell} \ZZ_\ell(n)$ commutes with limits. We have shown that $F_\ell(n)$ is indeed, as the notation suggests, the $n$-th $\ell$-adic Tate twist of $F_\ell:=F_\ell(0)$ when $\ell\neq p$. By definition, $F_\ell$ is $\ell$-adically complete, and $F/\ell$ is an $R/\ell$-constructible sheaf, so that $F_\ell(n)$ is an $R_{\ell}$-constructible sheaf (in the proétale sense as in \cite{Bhatt2015} or \cite{Hemo2023}).
\begin{prop}
	\begin{itemize}\label{prop:Weil_integral_model}
		\item[]
		\item 	Assume $\mathbf{L}_{F,n}$. Then $R\Gamma_{\mathrm{W,c}}(X,F(n))_{\ZZ_\ell} \simeq R\Gamma_{\mathrm{c}}(X_{\mathrm{proet}},F_\ell(n)).$
		\item Assume $\mathbf{R}(\dim X)$. The similar result holds for $R\Gamma_{\mathrm{Wh,c}}(X,F(n))$
	\end{itemize}
\end{prop}

\begin{proof}
	The functor $j_{!,\mathrm{proet}}$, for an open immersion $j$, has an exceptional left adjoint \cite[Cor. 6.1.5]{Bhatt2015} and thus commutes with limits, and derived global sections also commutes with limits. Since Weil-étale and Weil-eh cohomology with compact support agree with étale or proétale cohomology with compact support on torsion constructible sheaves, we obtain
	\[
	R\Gamma_{\mathrm{W,c}}(X,F(n))_{\ZZ_\ell}=R\Gamma_{\mathrm{W,c}}(X,F(n))^\wedge_\ell=\lim_k R\Gamma_{\mathrm{c}}(X_{\mathrm{proet}},F(n)/\ell^k)=R\Gamma_{\mathrm{c}}(X_{\mathrm{proet}},F_\ell(n))
	\]
	where the first equality follows from the perfectness over $\ZZ$ obtained from the assumption (and similarly for $R\Gamma_{\mathrm{Wh,c}}(X,F(n))$).
\end{proof}

\section{Reduction to the local case}\label{sec:reduction_local_case}

In this section, we explain how to reduce the proof of the main theorem \cref{thm:main} to a similar statement for the characteristic rational functions associated to the cohomology of $\ell$-adic sheaves.

We have a Verdier sequence 
\[\Perf(R)_{\mathrm{tor}}\to \Perf(R) \to \Perf(A),\]
hence $K(R,A)=K(\Perf(R)_{\mathrm{tor}})$ and similarly for $K_0$. The direct sum decomposition $\Perf(R)_{\mathrm{tor}}=\oplus_\ell \Perf(R)_{\ell^\infty\mathrm{-tor}}$ as well as the Verdier sequences $\Perf(R)_{\ell^\infty\mathrm{-tor}}\to \Perf(R_\ell) \to \Perf(A_\ell)$ now imply that $K(R,A)=\bigoplus_\ell K(R_\ell,A_\ell)$. Since \cref{thm:main} claims, for the special value, an equality in $K_0(R,A)=\bigoplus_\ell K_0(R_{\ell},A_{\ell})$, it suffices to prove the formula one prime at a time:
\begin{lem}\label{lem:reduction_local_case}
	For any prime $\ell$, let $i_\ell:\ct{A}^\times \hookrightarrow \ct{A_\ell}^\times$ denote the canonical inclusion and $\pi_\ell:K_0(R,A)\to K_0(R_{\ell},A_{\ell})$ denote the canonical projection. Assuming that $\chi_{\mathrm{W},X}(F(n))$ is well-defined, the special value formula is true if and only if 
	\[
	\delta_{\ell}(i_\ell(L^\ast_X(F,q^{-n})))=-\pi_\ell(\chi_{\mathrm{W},X}(F(n)))
	\]
	for all primes $\ell$. A similar statement holds for the Weil-eh case.
\end{lem}
We thus seek to reformulate the terms appearing in the above in terms intrinsic to the $\ell$-adic completion $F_\ell$ of $F$.

\subsection{Witte's noncommutative \texorpdfstring{$L$}{L}-functions}

\begin{lem}\label{lem:jabson_order}
	Let $B$ be a finite flat algebra over $\ZZ$. Then $J(B)$ is nilpotent; in particular $B$ is $J(B)$-adically complete.
\end{lem}
\begin{proof}
	Let $d$ be the rank of $B$, and let $\ell$ be a prime number. Then $B/\ell$ is a finite algebra of dimension $d$ over $\FF_\ell$, hence there is some integer $n\leq d$ such that $J(B/\ell)^n=J(B/\ell)^{n+1}$. By Nakayma's lemma we find $J(B/\ell)^n=0$ and a fortiori $J(B/\ell)^d=0$. Denote by $\pi_\ell: B\to B/\ell$ the canonical projection. We have $\pi_\ell(J(B)^d)\subseteq J(B/\ell)^d=0$, hence $J(B)^d\subseteq \bigcap_\ell (\ell)=0$ since $B$ is free over $\ZZ$.
\end{proof}
The above lemma implies that $\set{B[T]/(T^k)}$ and $\set{B[[T]]/J(B[[T]])^k}=\set{B[[T]]/(T,J(B))^k}$ are isomorphic pro-systems, so that $\lim K_1(B[T]/(T^k))\simeq \lim K_1(B[[T]]/J(B[[T]])^k)$.

\begin{lem}
	Let $B$ be a finite algebra over $\ZZ_p$. Then the $J(B)$-adic topology and the $p$-adic topology on $B$ are equivalent; in particular $B$ is an adic ring in the sense of \cite[Def. 3.1]{Witte2014}.
\end{lem}
\begin{proof}
	By \cite[5.9]{Lam2001}, we have $p \in J(B)$. Thus $J(B/p)=J(B)/p$ and the previous proof shows that $J(B)^d\subseteq (p)$ where $d$ is the dimension of $B/p$ over $\FF_p$. Finally, since $B$ is finite over $\ZZ_p$, it is $p$-adically complete.
\end{proof}

By the lemma, the ring $R_{\ell}$ is adic and $R_{\ell}[[T]]$ is again an adic ring whose Jacobson radical is $(T,J(R_\ell))$. Thus $K_1(R_{\ell}[[T]])\simeq\widehat{K}_1(R_\ell[[T]]):=\lim K_1(R_{\ell}[[T]]/(T,J(R_{\ell}))^k)$ \cite[Prop. 1.5.1]{Fukaya2006}. The $R_\ell$-constructible sheaf $F_\ell$ has an associated $L$-function element $\mathscr{L}_X(F_\ell)$ in $K_1(R_{\ell}[[T]])$ constructed as above through an Euler product \cite[Prop. 6.3]{Witte2014}. Moreover, if we denote by $L_X(F_\ell,T)$ the image under the reduced norm of $\mathscr{L}_X(F_\ell)\otimes_{R_\ell[[T]]} A_\ell[[T]]$ in $\ct{A_\ell}[[T]]^\times$, then 
\begin{equation}\label{eq:L_function_equals_l_adic_L_function}
	i_\ell(L_X(F,T))=L_X(F_\ell,T).
\end{equation} 
where $i_\ell:\ct{A}^\times \hookrightarrow \ct{A_\ell}^\times$ denotes the canonical inclusion. Indeed, the pullback to a closed point $x$ in the proétale topology commutes with limits \cite[Cor. 6.1.5]{Bhatt2015}, hence $(F_\ell)_x=(F_x)^\wedge_\ell=F_x\otimes_\ZZ \ZZ_\ell$, the latter equality by the perfectness of $F_x$.

\subsection{Comparison of the integral and \texorpdfstring{$\ell$}{l}-adic structures}

There is again an element $e_\ell:=e\otimes 1\in H^1_W(\Spec(k),\ZZ)_{\ZZ_\ell}\simeq H^1(\Spec(k),\ZZ_\ell)\simeq \Hom_{\mathrm{cont}}(\Gal(k^{\mathrm{sep}}/k),\ZZ_\ell)$, described by $\varphi\mapsto 1$ and corresponding to the extension of discrete $\Gal(k^{\mathrm{sep}}/k)$-modules
\[
\ZZ_\ell \to \ZZ_\ell^2 \to \ZZ_\ell
\]
with the similar action as before.

\begin{defi}\label{def:local_Euler_char}
	Define a trivialization of $[R\Gamma_{\mathrm{c}}(X,F_\ell(n))_{\QQ_\ell}]_{A_\ell}$ by
	\[
	\alpha_\ell:[R\Gamma_{\mathrm{c}}(X,F_\ell(n))_{\QQ_\ell}]_{A_\ell}\xrightarrow{\tau^\geq} \bigotimes_{i\in \ZZ}[H^i_{\mathrm{c}}(X,F_\ell(n))_{\QQ_\ell}]_{A_\ell}^{(-1)^i} \xleftarrow[\simeq]{\sigma^{\geq}} [(H^\ast_{\mathrm{c}}(X,F_\ell(n))_{\QQ_\ell},\cup e_\ell)]_{A_\ell} \xrightarrow{[0]} \mathbf{1}_{A_\ell},
	\]
	where the map $\tau^\geq$ is obtained from the canonical truncation filtration, the map $\sigma^\geq$ from the use of the stupid truncation filtration, and the map $[0]$ is the isomorphism to $\mathbf{1}$ for an acyclic object or complex.
	
	This defines an Euler characteristic $\chi_X(F_\ell(n)):=\chi([R\Gamma_{\mathrm{c}}(X,F_\ell(n))]_{R_\ell},\alpha_\ell)\in K_0(R_\ell,A_\ell)$.
\end{defi}

\begin{lem}\label{two_trivializations}
	\begin{itemize}
		\item[]
		\item Assume $\mathbf{L}_{F,n}$. Then we have $\alpha_\ell=\alpha^\mathrm{W}\otimes 1$.
		\item Assume $\mathbf{R}(\dim X)$. Then we have $\alpha_\ell=\alpha^\mathrm{Wh}\otimes 1$.
	\end{itemize}

\end{lem}
\begin{proof}
	We treat the first case. Under $\mathbf{L}_{F,n}$, we can consider the following diagram.
\[\adjustbox{center}{	\resizebox{\textwidth}{!}{
		\begin{tikzcd}[ampersand replacement=\&,column sep=small,cramped]
		{[R\Gamma_{\mathrm{W,c}}(X,F(n))_\QQ]} \& {\bigotimes_{i\in \ZZ}[H^i_{\mathrm{W,c}}(X,F(n))_\QQ]^{(-1)^i}} \& {[(H^\ast_{\mathrm{W,c}}(X,F(n))_\QQ,\cup e)]} \\
		{[V[-1]\oplus V ]} \& {\bigotimes_{i\in \ZZ}[H^{i-1}V\oplus H^iV]^{(-1)^i}} \& \begin{array}{c} [(H^{\ast-1}V \oplus H^\ast V,\begin{pmatrix} 0 & 1\\ 0 & 0\end{pmatrix})]  \end{array} \\
		{[V]^{-1}[V]} \& {\bigotimes_{i\in \ZZ}[H^{i}V]^{(-1)^{i+1}}} \otimes \bigotimes_{i\in \ZZ}[H^iV]^{(-1)^i} \& {[(H^{\ast-1},0)][(H^{\ast},0)]} \\
		\\
		\&\&\&\& {\mathbf{1}}
		\arrow["{\tau^\geq}", from=1-1, to=1-2]
		\arrow[from=1-1, to=2-1]
		\arrow[from=1-2, to=2-2]
		\arrow["{\sigma^\geq}"', from=1-3, to=1-2]
		\arrow[from=1-3, to=2-3]
		\arrow["{[0]}", from=1-3, to=5-5,bend left=30]
		\arrow["{\tau^\geq}", from=2-1, to=2-2]
		\arrow[from=2-1, to=3-1]
		\arrow["\ast"', from=2-2, to=3-2]
		\arrow["{\sigma^\geq}"', from=2-3, to=2-2]
		\arrow["{\mathrm{cone(id)}}"', from=2-3, to=3-3]
		\arrow["{[0]}"{description}, from=2-3, to=5-5,bend left=15]
		\arrow["{\tau^\geq}", from=3-1, to=3-2]
		\arrow["{\sigma^\geq}"', from=3-3, to=3-2]
		\arrow["{\mathrm{can}}"', from=3-1, to=5-5, bend right=8]
		\arrow["{\mathrm{mult}\circ\otimes_i\mathrm{can}^{(-1)^i}}"{description, pos=0.3},"\ast", from=3-2, to=5-5]
		\arrow["{\mathrm{can}}"{description}, from=3-3, to=5-5]
		\end{tikzcd}
}}\]

Starred arrows denote an implicit use of the commutativity constrain. The map $\mathrm{can}$ for an object $X$ is the one witnessing that $X^{-1}$ is the inverse of $X$; if $X=[M]$, it is obtained from the fiber sequence $M\to 0 \to M[1]$ under the identification $[M[1]]=[M]^{-1}$. The map denoted $\mathrm{cone(id)}$ is obtained from the levelwise split short exact sequence of complexes
\[
0\to (H^{\ast-1}V, 0) \to (H^{\ast-1}V \oplus H^\ast V,\begin{pmatrix} 0 & 1\\ 0 & 0\end{pmatrix}) \to (H^\ast V,0) \to 0
\]
realizing the cone construction on $\mathrm{id}: (H^{\ast-1}V, 0) \to  (H^{\ast-1}V, 0)$.

The diagram commutes, hence the trivialization $\alpha^\mathrm{W}$ of $R\Gamma_{\mathrm{W,c}}(X,F(n))_\QQ$ is also described by the top line postcomposed with $[0]$. The result follows since Weil-étale cohomology is an integral model of $\ell$-adic cohomology (\cref{prop:Weil_integral_model}).
\end{proof}

\begin{prop}\label{prop:reduction_local_case}
	\begin{itemize}
		\item[]
		\item Assume $\mathbf{L}_{F,n}$.  Then
		\[
		\pi_\ell(\chi_{\mathrm{W},X}(F(n)))=\chi_X(F_\ell(n)).
		\]
		\item Assume $\mathbf{R}(\dim X)$.  Then
		\[
		\pi_\ell(\chi_{\mathrm{Wh},X}(F(n)))=\chi_X(F_\ell(n)).
		\]
	\end{itemize}
\end{prop}

\begin{proof}
	This follows from the previous lemma by the compatibility of the formation of refined Euler characteristics of trivialized complexes with base change.
\end{proof}

\subsection{Special values of \texorpdfstring{$\ell$}{l}-adic power series and rationality}

\begin{defi}
	Let $K$ be a finite extension of $\QQ_\ell$ with ring of integers $\cal{O}_K$, and let $B$ be an $\cal{O}_K$-order in a semi-simple algebra $B_K$ over $K$. Let $\eta\in \ct{B_K}[[T]]$, let $x\in \cal{O}_K\backslash\{0\}$, and write $\ct{B_K}\simeq \prod_{i=1}^n K_i$ where each $K_i$ is an $\ell$-adic local field. We say that $\eta$ has a well-defined special value at $t=x$ if each component $\eta_i$ of $\eta$ in $K_i[[T]]$ defines an $\ell$-adic meromorphic function $\eta_i$ on the closed unit disk without pole at $0$, i.e. if $\eta_i\in K_i\langle T\rangle_{(T)}\subset K_i[[T]]$ where $K_i\langle T\rangle$ is the Tate algebra in one variable over $K_i$, consisting of power series whose sequence of coefficients tends to $0$. Then by the Weierstrass preparation theorem, one can write
		\[
		\eta_i=(x-T)^{r_i}\frac{P_i}{Q_i}u_i,
		\]
		with $P_i,Q_i\in K_i[T]$ prime to $(x-T)$ and $u_i\in K_i\langle T\rangle^\times$.
	The vanishing order of $\eta$ at $t=x$ is then the element $\ord_{t=x}\eta:=(r_i)\in \ZZ^n=K_0(A_\ell)$ and the special value $\eta^\ast(x)\in \ct{A_\ell}^\times$ of $\eta$ at $t=x$ is the element $(\frac{P_i(x)}{Q_i(x)}u_i(x))$.\footnote{Again, we have chosen the normalization where the leading term has the form $\eta^\ast(x)(x-T)^{\ord_{t=x}\eta}$.}
	
%
\end{defi}
\begin{rmk}
	\begin{itemize}
		\item[]
		\item This definition is not the most general but it will be enough for our purposes.
		\item The product of localizations $\prod_i K_i\langle T\rangle_{(T)}$ is the localization, denoted $\ct{B_K}\langle T\rangle)$, of $\ct{B_K}\langle T\rangle$ at elements that become invertible when evaluated at $0$. The definition thus asks that $\eta\in L_0(\ct{B_K}\langle T\rangle)$; then one can write $\eta=(x-T)^{\underline{\mathbf{r}}}\cdot\frac{P}{Q}\cdot u,$ with $P,Q\in \ct{B_K}[T]$ prime\footnote{On each component.} to $(x-T)$, $u\in \ct{B_K}\langle T\rangle^\times$ and $\underline{\mathbf{r}}\in K_0(B_K)\xrightarrow[\mathrm{rrank}_{B_K}]{\simeq}\ZZ^n$.
			
			With this formulation we have $\ord_{t=x}\eta=\underline{\mathbf{r}}$ and $\eta^\ast(x)=\frac{P(x)}{Q(x)}u(x)$.
	\end{itemize}
\end{rmk}

Recall that for a prime $\ell$, we let $i_\ell:\ct{A}^\times \hookrightarrow \ct{A_\ell}^\times$ denote the canonical inclusion.
\begin{lem}\label{lem:rationality}
	Let $\eta\in \ct{A}[[T]]^\times$ and assume that there exists a prime $\ell$ such that $i_\ell(\eta)\in \ct{A_\ell}[[T]]^\times$ is rational (that is, every component in the decomposition of $\ct{A_\ell}$ into a product of local fields is a rational function). Then:
	\begin{enumerate}
		\item $\eta$ is rational and its special value at $t=x$ is well-defined for all $x\in \QQ\backslash\set{0}$;
		\item for every prime $\ell'$, $i_{\ell'}(\eta)\in \ct{A_{\ell'}}[[T]]^\times$ is rational, its special value at $t=x$ is well-defined for all $x\in \QQ\backslash\set{0}$ and equals the image of $\eta^\ast(x)$ inside $\ct{A_{\ell'}}^\times$, and its vanishing order at $t=q^{-n}$ equals the image of $\ord_{t=x}\eta$ in $K_0(A_{\ell'})$.
	\end{enumerate} 
\end{lem}
\begin{proof}
	Everything follows readily once it has been observed that $\eta$ is rational, because the vanishing order and special value of a rational function are completely algebraic. Write the centers $\ct{A}=\prod_i K_i$ and $\ct{A_\ell}=\prod_i \prod_j L_{i,j}$ as products of fields, such that each $L_{i,j}$ is an extension of $K_i$. Each component $\eta_i$ is rational if and only if the adequate Hankel determinants vanish, which can be detected after base change to some $L_{i,j}$; but this is then the Hankel determinants for the component $i_{\ell}(\eta)_{i,j}$ of $i_{\ell}(\eta)$, which is eventually zero (for a fixed, big enough size) by the assumption that $\eta_\ell$ is rational.
\end{proof}

\subsection{Witte's theorem on the rationality of the \texorpdfstring{$L$}{L}-function for \texorpdfstring{$\ell\neq p$}{l!=p}}\label{subsec:rationality}
In this section, we use Witte's theorem on the rationality of non-commutative $\ell$-adic $L$-functions for $\ell\neq p$ \cite{Witte2014} to reduce the local problem to the computation of the special value of a characteristic rational function, namely that associated to the $\ell$-adic cohomology with compact support of the considered coefficients.

Let $\ell$ be a prime, let $K$ be a finite extension of $\QQ_\ell$ with ring of integers $\cal{O}_K$, and let $B$ be an $\cal{O}_K$-order in a semi-simple algebra $B_K$ over $K$. Let $G$ be a $B$-constructible sheaf on $X$ in the proétale sense as in \cite{Bhatt2015} or \cite{Hemo2023}. The geometric Frobenius $\phi$ acts on $R\Gamma_{\mathrm{c}}(X_{\overline{k}},G)$ and $1-\phi\cdot T$ is an endomorphism of $R\Gamma_{\mathrm{c}}(X_{\overline{k}},G)\otimes_{B} B[T]$ that becomes an isomorphism when evaluated at $0$. For an associative ring $B$, let us denote by $L_0 B[T]$ the Cohn localization of $B[T]$ at matrices that become invertible when evaluated at $0$. We get that $1-\phi\cdot T$ is an automorphism of $R\Gamma_{\mathrm{c}}(X_{\overline{k}},G)\otimes_{B} L_0 B[T]\in \Perf(L_0 B[T])$, and thus defines a class $\langle 1-\phi\cdot T \mid R\Gamma_{\mathrm{c}}(X_{\overline{k}},G)\otimes_{B} L_0 B[T]\rangle \in K_1(L_0 B[T])$. For commutative rings or semi-simple algebras, the invertibility of a matrix can be checked on its determinant or reduced norm, so the Cohn localization is the usual localization inverting elements in the ring, and it is an Ore localization.

In particular, if $\ct{B_K}\simeq \prod_{i=1}^d K_i$ is a decomposition into a product of fields, then $L_0(\ct{B_K}[T])=\prod_{i=1}^d K_i[T]_{(T)}\subseteq \ct{B_K}[[T]]$ is the ring of power series that are rational without pole at $0$ on each component. Moreover, by construction, if $C=M_n(B)$ is a ring of matrices over an associative ring $B$, then $L_0 M_n(B)=M_n(L_0B)$. Applying this to the components of $B_K\otimes_{\QQ_{\ell}} \overline{{\QQ_{\ell}}}\simeq \prod_{i=1}^m M_{k_i}(\overline{{\QQ_{\ell}}})$, we get by Morita invariance a map
\begin{align*}
K_1(L_0B_K[T])\to K_1\left(L_0 ((B_K\otimes_{\QQ_{\ell}} \overline{{\QQ_{\ell}}})[T])\right)=\prod_{i=1}^m K_1( \overline{{\QQ_{\ell}}}[T]_{(T)})&\xrightarrow[\prod \det]{\simeq}\prod_{i=1}^m(\overline{{\QQ_{\ell}}}[T]_{(T)})^\times\\&=\left(L_0(\ct{B_K\otimes_{\QQ_{\ell}} \overline{{\QQ_{\ell}}})[T]}\right)^\times
\end{align*}
We claim that this map factors as
\[\begin{tikzcd}[ampersand replacement=\&]
	{K_1(L_0B_K[T])} \& {K_1\left(L_0 ((B_K\otimes_{\QQ_{\ell}} \overline{{\QQ_{\ell}}})[T])\right)} \\
	{\bigl(L_0(\ct{B_K}[T])\bigr)^\times} \& {\left(L_0(\ct{B_K\otimes_{\QQ_{\ell}} \overline{{\QQ_{\ell}}})[T]}\right)^\times}
	\arrow[from=1-1, to=1-2]
	\arrow["{\mathrm{Nrd}_{B_K}}"', dashed, from=1-1, to=2-1]
	\arrow[from=1-2, to=2-2]
	\arrow[hook, from=2-1, to=2-2]
\end{tikzcd}\]
giving a \emph{reduced norm} map $\mathrm{Nrd}_{B_K}:K_1(L_0B_K[T])\to \left(L_0(\ct{B_K}[T])\right)^\times$. Indeed, it suffices to show that Galois-invariant rational functions over $\ct{B_K\otimes_{\QQ_{\ell}} \overline{{\QQ_{\ell}}}}$ are rational functions over $\ct{B_K}$. This is known for power series since it can be checked on coefficients, and then it suffices to show that a power series over $\ct{B_K}$ which is rational over $\ct{B_K\otimes_{\QQ_{\ell}} \overline{{\QQ_{\ell}}}}$ is rational; this is the Hankel determinant argument appearing in the proof of \cref{lem:rationality}.

\begin{defi}
	Let $G$ be a $B$-constructible sheaf on $X$. The cohomological $L$-function of $G$ is 
	\[
	\xi_X(G,T):=\mathrm{Nrd}_{B_K}(1- \phi \cdot T \mid R\Gamma_{\mathrm{c}}(X_{\overline{k}},G)\otimes_{B} L_0(B_K[T]))^{-1}\in (L_0\ct{B_K}[T])^\times.
	\]
	where $\phi$ denotes the geometric Frobenius. It is a rational function without a pole or zero at $0$.
\end{defi}

\begin{thm}[\cite{sga412fonctionsmodl},{\cite[Thm. 1.1]{Witte2014}}]\label{thm:rationality}
	Suppose $\ell\neq p$. The $L$-function element $\mathscr{L}_X(G)$ has the following cohomological interpretation:
	\[
	\mathscr{L}_X(G)=\langle 1- \phi \cdot T \mid R\Gamma_{\mathrm{c}}(X_{\overline{k}},G)\otimes_{B} L_0(B[T])\rangle ^{-1}\otimes_{L_0(B[T])} B[[T]]
	\]
	In particular, $L_X(G,T)=\xi_X(G,T)$ is rational without a zero or pole at $0$.
\end{thm}
\begin{proof}
We want to aply \cite[Thm. 1.1]{Witte2014}. We claim an equality in 
\[
K_1(B[[T]])=\lim K_1(B[[T]]/(J(B),T)^n)
\]
\cite[Prop. 1.5.1]{Fukaya2006}. There is an integer $k\geq 1$ such that $\ell^k\in J(B)$, and thus comparison maps 
\[
K_1(B[[T]])\to K_1(R/\ell^k[[T]])\to K_1(B[[T]]/(J(B),T)^n)
\]
We have $G/\ell^k\in\Cons(X,R/\ell^k)$ \cite[Lem. 6.3.12]{Bhatt2015}, and the images of the two elements of the claim in $K_1(R/\ell^k[[T]])$ only depend on the class $[G/\ell^k]\in K_0(\Cons(X,R/\ell^k))$, hence it suffices to show the claim for $[G/\ell^k]$. Now observe that any constructible complex over $R/\ell^k$ has finite-type--constructible cohomology sheaves because $R/\ell^k$ is Noetherian, and has finite Tor dimension as this can be checked on each strata (by the projection formula) and locally. Moreover, any such complex is represented by a finite complex of flat constructible sheaves by \cite[Prop.-def. 4.6]{sga412rapport}. Conversely, any such complex is clearly constructible. 

Let $\mathbf{PDG}(X,R/\ell^k)$ denote Witte's Waldhausen category of DG-flat constructible complexes. Given a level-wise injective map of complexes, the corresponding cofibration sequence is a short exact sequence of complexes hence is isomorphic to a cofiber sequence in the derived $\infty$-category. This gives a well-defined map $K_0(\mathbf{PDG}(X,R/\ell^k))\to K_0(\Cons(X,R/\ell^k))$. A finite complex of flat sheaves is DG-flat, so the map is surjective by the previous discussion. Since the maps
\[\begin{array}{cl}
	K_0(\mathbf{PDG}(X,R/\ell^k)) &\to K_1(R/\ell^k[[T]])\\
	{[G]} &\mapsto \prod_{x\in X_0} \langle 1-\phi\cdot T \mid G_x\otimes_{R/\ell^k}R/\ell^k[[T]] \rangle^{-1},
\end{array}\]

\[\begin{array}{cl}
	K_0(\mathbf{PDG}(X,R/\ell^k)) &\to K_1(R/\ell^k[[T]])\\
	{[G]} &\mapsto \langle 1-\phi\cdot T \mid R\Gamma_c(X_{\overline{k}},G)\otimes_{R/\ell^k}R/\ell^k[[T]] \rangle^{-1},
\end{array}\]
factor through $K_0(\Cons(X,R/\ell^k))$, we conclude by applying \cite[Thm. 1.1]{Witte2014} to a lift of $[G/\ell^k]$.

\end{proof}

Assume $\mathbf{L}_{F,n}$. Then $(R\Gamma_{\mathrm{W,c}}(X,F)_\QQ)_{\QQ_\ell} = R\Gamma_{\mathrm{c}}(X,F_\ell(n))_{\QQ_\ell}$ so by compatibility of the formation of ranks with field extensions, we get $\mathrm{rrank}_{A} H^i_{\mathrm{W,c}}(X,F(n))_{\QQ} = \mathrm{rrank}_{A_\ell} H^i_{\mathrm{c}}(X,F_\ell)_{\QQ_\ell}$ for all $i\in \ZZ$. The similar argument applies for Weil-eh cohomology, assuming $\mathbf{R}(\dim X)$. We can thus deduce the following reduction to the $\ell$-adic case:

\begin{cor}\label{cor:reduction_local_case_0}
	Let $F$ be an $R$-constructible complex on a finite type separated scheme $X$ over $k$.
	\begin{enumerate}
		\item The $L$-function $L_X(F,T)$ is rational and has a well-defined special value at $t=x$ for all $x\in \QQ\backslash\set{0}$. We have
		\[
		\ord_{t=q^{-n}}L_X(F,T)=\ord_{t=q^{-n}} \xi_X(F_\ell,T)
		\]
		for any prime $\ell\neq p$. 
		\item Let $n\leq 0$ and assume $\mathbf{L}_{F,n}$. Then the vanishing order formula at $t=q^{-n}$ for $L_X(F,T)$ is true if and only if we have the vanishing formula for $\xi_X(F_\ell,T)$:
		\[
		\ord_{t=q^{-n}} \xi_X(F_\ell,T)= \sum_i (-1)^i\cdot i \cdot \mathrm{rrank}_{A_\ell} H^i_{\mathrm{c}}(X,F_\ell(n))_{\QQ_\ell} 
		\]
		for some $\ell\neq p$
		\item The special value formula at $t=q^{-n}$ for $L_X(F,T)$ is true if and only if we have the special value formula for $\xi_X(F_\ell,T)$:
		\[
		\delta_\ell(\xi_X^\ast(F_\ell,q^{-n}))=-\chi_X(F_\ell(n))
		\]
		for $\ell\neq p$ and
		\[
		\delta_p(L^\ast(F_p,q^{-n}))=-\chi_X(F_p(n)).
		\]
		\item Assume $\mathbf{R}(\dim X)$. Then the similar result holds using $R\Gamma_{\mathrm{Wh,c}}(X,F(n))$.
	\end{enumerate}
\end{cor}
\begin{proof}
	We have $i_\ell(L_X(F,T))=L_X(F_\ell,T)$ (\cref{eq:L_function_equals_l_adic_L_function}) so \cref{lem:rationality} applied to \cref{thm:rationality} gives the rationality and the well-definedness of special values, and the identification of the vanishing order with that of the cohomological $L$-function for any prime $\ell\neq p$. Assume $\mathbf{L}_{F,n}$.  We have $R\Gamma_{\mathrm{W,c}}(X,F(n))\in\Perf(R)$ by \cref{cor:perfect_Z_implies_perfect_R}. Thus, the Weil-étale Euler characteristic $\chi_{\mathrm{W},X}(F(n))$ of \cref{def:weil_etale_euler_char} is well-defined. The statement for the vanishing order formula then follows from the previous discussion and point 1., and the one for the special value is \cref{lem:reduction_local_case,prop:reduction_local_case}. The alternative statement for the Weil-éh cohomology is proven similarly.
\end{proof}

\subsection{The case \texorpdfstring{$\ell=p$}{l=p}}\label{subsec:l=p}

In this section, we complete the reduction to the computation of the special value of the cohomological $L$-function, in the case $\ell=p$. Let $K$ be a finite extension of $\QQ_p$ with ring of integers $\cal{O}_K$, and let $B$ be an $\cal{O}_K$-order in a semi-simple algebra $B_K$ over $K$. 
For a $B$-constructible sheaf $G$ on $X$, the $L$-function $L_X(G,T)$ is not equal anymore to $\xi_X(G,T)$. Witte's results in \cite{Witte2016}, generalizing \cite{Emerton2001} to the non-commutative case, quantifies their difference; in particular, after applying the boundary $\delta: \ct{B_K}^\times \to K_0(B,B_K)$, their special value at $q^{-n}$ coincide when $n\leq 0$:

\begin{thm}\label{thm:same_vs_under_boundary}
	Let $G$ be a $B$-constructible complex on a finite type separated scheme $X$ over $k$ and let $x\in \cal{O}_K\backslash\{0\}$. We have
	\[
	\delta(L_X^\ast(G,x))=\delta(\xi_X^\ast(G,x)) \in K_0(B,B_K).
	\]
	This applies in particular for $x=q^{-n}$ when $n\leq 0$.
\end{thm}

Consider the Tate algebra $B\langle T\rangle:= \lim B/p^k B [T]$, which is identified with the sub-ring of $B[[T]]$ of power series with coefficients converging to $0$ in the $p$-adic topology. We can also define the Tate algebras $B_K\langle T\rangle$ and $\ct{B_K}\langle T \rangle$ as the sub-algebras of the respective formal power series algebras where coefficients converge to $0$ in the $p$-adic topology.

The evaluation at any $x\in \cal{O}_K$ is a well-defined continuous map  $\ct{B_K}\langle T \rangle \xrightarrow{\mathrm{ev}_x} \ct{B_K}$. In particular, if $\eta\in \ct{B_K}\langle T \rangle^\times$, then $\eta$ has no zero or poles on the closed unit disk, so that $\eta$ has a well-defined special value at $x$ given by its evaluation at $x$, and vanishing order $0$. Moreover:

\begin{lem}
	Let $\eta\in K_1(B\langle T\rangle)$. Then $\mathrm{Nrd}_{B_K}(\eta \otimes_{B\langle T\rangle} B_K\langle T\rangle)\in \ct{B_K}\langle T \rangle^\times$ has a well-defined special value at $x$, with order of vanishing $0$ and special value $\eta^\ast(x)\in \ker(\delta)$.
\end{lem}

\begin{proof}
	The previous discussion shows that the special value is given by evaluating at $x$, and the vanishing order is $0$.
	The following diagram commutes
	\[\begin{tikzcd}[ampersand replacement=\&]
		{K_1(B\langle T\rangle)} \& {K_1(B_K\langle T\rangle)} \& {\ct{B_K}\langle T\rangle^\times} \\
		{K_1(B)} \& {K_1(B_K)} \& {\ct{B_K}^\times}
		\arrow[from=1-1, to=1-2]
		\arrow["{\mathrm{ev}_x}", from=1-1, to=2-1]
		\arrow["{\mathrm{Nrd}}", from=1-2, to=1-3]
		\arrow["{\mathrm{ev}_x}", from=1-2, to=2-2]
		\arrow["{\mathrm{ev}_x}", from=1-3, to=2-3]
		\arrow[from=2-1, to=2-2]
		\arrow["{\mathrm{Nrd}}"', "\simeq", from=2-2, to=2-3]
	\end{tikzcd}\]
	because after base change to $B_K\otimes_{\QQ_p}\overline{\QQ_p}$, the reduced norm is a determinant.
	We thus have $\left(\mathrm{Nrd}_{B_K}(\eta\otimes_{B\langle T \rangle} 1)\right)^\ast(x)=\mathrm{ev}_{x}\left(\mathrm{Nrd}_{B_K}(\eta\otimes_{B\langle T \rangle} 1)\right)=\mathrm{Nrd}_{B_K}(\mathrm{ev}_x(\eta)\otimes_{B} B_K)$. Since $\ker(\delta)=\mathrm{Nrd}_{B_K}(K_1(B))$, we deduce the result.
\end{proof}

\begin{thm}[{\cite[Thm. 1.1, Prop. 2.2]{Witte2016}}]\label{thm:witte2}
	Let $G$ be a $B$-constructible complex on a finite type separated scheme $X$ over $k$. Then the quotient 
	\[
	\frac{\mathscr{L}_X(G)}{\langle 1-\phi\cdot T\mid R\Gamma_{\mathrm{c}}(X_{\overline{k}},G) \otimes_{B} B[[T]] \rangle^{-1}}\in K_1(B[[T]])
	\]
	lifts to $K_1(B\langle T\rangle)$.
\end{thm}
\begin{proof}
	The proof is the same as for \cref{thm:rationality}, noting that there is a surjection from $K_1(B\langle T\rangle)$ to $\widehat{K}_1(B\langle T\rangle):=\lim K_1(B/J(B)^n[T])$ \cite[Prop. 2.2]{Witte2016}.
\end{proof}

We can now apply the lemma to prove \cref{thm:same_vs_under_boundary}:
\begin{proof}[Proof of \cref{thm:same_vs_under_boundary}]
	We already know that $L_X(G,T)$ and $\xi_X(G,T)$ are rational, hence have a well-defined special value at $x$. It is clear from our definitions that special values are multiplicative, so by the lemma applied to a lift in $K_1(B\langle T\rangle)$ of $L_X(G,T)/\xi_X(G,T)$, we obtain
	\[
	\delta(L_X^\ast(G,x))-\delta(\xi_X^\ast(G,x))=\delta\Bigl(\bigl(\frac{L_X(G,T)}{\xi_X(G,T)}\bigr)^\ast(x)\Bigr)=0 \qedhere
	\]
\end{proof}

We can now generalize \cref{cor:reduction_local_case_0} to include the case $\ell=p$:

\begin{cor}\label{cor:reduction_local_case}
	Let $F$ be an $R$-constructible complex on a finite type separated scheme $X$ over $k$.
	\begin{enumerate}
		\item The $L$-function $L_X(F,T)$ is rational and has a well-defined special value at $t=x$ for all $x\in \QQ\backslash\set{0}$. We have
		\[
		\ord_{t=q^{-n}}L_X(F,T)=\ord_{t=q^{-n}} \xi_X(F_\ell,T)
		\]
		for any prime $\ell\neq p$. 
		\item Let $n\leq 0$ and assume $\mathbf{L}_{F,n}$. Then the vanishing order formula at $t=q^{-n}$ for $L_X(F,T)$ is true if and only if we have the vanishing formula for $\xi_X(F_\ell,T)$:
		\[
		\ord_{t=q^{-n}} \xi_X(F_\ell,T)= \sum_i (-1)^i\cdot i \cdot \mathrm{rrank}_{A_\ell} H^i_{\mathrm{c}}(X,F_\ell(n))_{\QQ_\ell} 
		\]
		for some $\ell\neq p$
		\item The special value formula at $t=q^{-n}$ for $L_X(F,T)$ is true if and only if we have the special value formula for $\xi_X(F_\ell,T)$:
		\[
		\delta_\ell(\xi_X^\ast(F_\ell,q^{-n}))=-\chi_X(F_\ell(n))
		\]
		for \emph{all} primes $\ell$.
		\item Assume $\mathbf{R}(\dim X)$. Then the similar result holds using $R\Gamma_{\mathrm{Wh,c}}(X,F(n))$.
	\end{enumerate}
\end{cor}
\begin{proof}
	Combine \ref{cor:reduction_local_case_0} and \cref{thm:same_vs_under_boundary}.
\end{proof}

\section{Semi-simplicity}\label{sec:semisimplicity}

In this section, we recall the definition of semi-simplicity at $0$ of an endomorphism, and provide equivalent characterizations. We then use this to show that the Weil-étale cohomology forces semi-simplicity at $0$ of $1-\phi$ on $\ell$-adic cohomology when it is an integral model of the latter, i.e. when it is perfect.

\begin{defi}
	Let $B$ be an associative unital ring. We say that a pair $(M,\theta)$ of a left $B$-module $M$ with endomorphism $\theta$ (resp. a complex $M\in D(B)$ with endomorphism $\theta$) is semi-simple at $0$ if $M$ admits a $\theta$-equivariant direct sum decomposition $M=V\oplus W$ in $B$-modules (resp. in $D(B)$) in which $\theta_{|V}=0$ and $\theta_{|W}$ is an automorphism of $W$.
\end{defi} 
If a complex $M\in D(B)$ is semi-simple at $0$, the fiber of $\theta$ inherits via the decomposition a direct sum structure $\fib(\theta)=V[-1]\oplus V$.

There is a fiber sequence
\[
R\Gamma_{\mathrm{c}}(X,F_\ell(n))\to R\Gamma_{\mathrm{c}}(X_{\overline{k}},F_\ell(n)) \xrightarrow{1-\phi} R\Gamma_{\mathrm{c}}(X_{\overline{k}},F_\ell(n))
\]
involving the \emph{geometric} Frobenius $\phi$, inducing a map
\[
\beta_0:R\Gamma_{\mathrm{c}}(X,F_\ell(n)) \to R\Gamma_{\mathrm{c}}(X_{\overline{k}},F_\ell(n)) \to R\Gamma_{\mathrm{c}}(X,F_\ell(n))[1]
\]

We will use the following:
\begin{lem}[{\cite[Prop. 3.2 (i)]{Burns}}]\label{lem:semisimplicity}
	Let $B$ be a semi-simple algebra over a field and let $M$ be a finitely generated left $B$-module with an endomorphism $\theta$. The following are equivalent: 
	\begin{enumerate}
		\item $\theta$ is semi-simple at $0$;
		\item the canonical map $\ker(\theta)\hookrightarrow M \twoheadrightarrow \mathrm{coker}(\theta)$ is an isomorphism.
	\end{enumerate}
	
	Let $M\in \Perf(B)$ be a perfect complex with an endomorphism $\theta$ and let $F:=\fib(\theta)$. Then we have a natural map $\beta_0:F\to M \to F[1]$ such that $(H^\ast F,H^\ast \beta)$ is a complex. The following are equivalent:
	\begin{enumerate}
		\item $\theta$ is semi-simple at $0$;
		\item $H^i(\theta)$ is semi-simple at $0$ for all $i\in \ZZ$;
		\item The complex $(H^\ast F,H^\ast \beta_0)$ is acyclic.
	\end{enumerate}
\end{lem}

\begin{proof}
	The first claim is standard linear algebra: since $B$ is semi-simple, we can choose a direct summand for $\ker\theta$ in $M$, and 2. implies both that this summand is actually a $\theta$-equivariant summand, and that the restriction of $\theta$ to that summand is an isomorphism.
	
	Let us prove the second claim. Since $B$ is semi-simple, we have $M\simeq \oplus_i H^i(M)[-i]$ and under that isomorphism, $\theta$ corresponds necessarily to $\oplus H^i(\theta)[-i]$ because there are no $\Ext^1$. This shows the equivalence between 1. and 2. Let $V^i:=\ker H^i\theta$, $W^i:=\mathrm{coker} H^i\theta$, and denote by $\mathrm{can}^i$ the composite $V^i\to H^iM\to W^i$. From the long exact sequence associated to the fiber sequence $F\to M\to M$, we deduce that the following diagram, where the diagonal maps denote the mono-epi factorizations of the middle row, commutes:
	\[\begin{tikzcd}[ampersand replacement=\&,sep=small]
		{H^{i-1}F} \&\& {H^{i-1}M} \&\& {H^iF} \&\& {H^iM} \&\& {H^{i+1}F} \\
		\& {V^{i-1}} \&\& {W^{i-1}} \&\& {V^{i}} \&\& {W^{i}}
		\arrow[from=1-1, to=1-3]
		\arrow["{H^{i-1}\beta_0}", bend left= 15, from=1-1, to=1-5]
		\arrow[two heads, from=1-1, to=2-2]
		\arrow[from=1-3, to=1-5]
		\arrow[two heads, from=1-3, to=2-4]
		\arrow[from=1-5, to=1-7]
		\arrow["{H^i\beta_0}", bend left=15, from=1-5, to=1-9]
		\arrow[two heads, "v^i", from=1-5, to=2-6]
		\arrow[from=1-7, to=1-9]
		\arrow[two heads, from=1-7, to=2-8]
		\arrow[hook, from=2-2, to=1-3]
		\arrow["{\mathrm{can}^{i-1}}", from=2-2, to=2-4]
		\arrow[hook, "u^i", from=2-4, to=1-5]
		\arrow[hook, from=2-6, to=1-7]
		\arrow["{\mathrm{can}^i}", from=2-6, to=2-8]
		\arrow[hook, from=2-8, to=1-9]
	\end{tikzcd}\]
		Since $H^{i-1}M \to H^iF\to H^iM$ is exact, we have $\ker (v^i)= u^i(W^{i-1})\simeq W^{i-1}$. The statement we are after only depends on the category of modules of $B$ so we can assume that $B$ is a division algebra and thus reason with dimensions. Put $k_i=\dim V^i$, $\ell_i=\dim W^i$, $\alpha_i=\dim(\ker(\mathrm{can}^i))$, $\beta_i=\rank(\mathrm{can}^i)$. Then on the one hand  $\beta_{i-1}\leq  \ell_{i-1}$, and on the other hand we compute by inspecting the diagram that $\rank (H^{i-1}\beta_0)=\beta_{i-1}$ and $\dim(\ker H^i\beta_0)=\alpha_i+\ell_{i-1}$. From this, we deduce that $(H^\ast F,H^\ast \beta_0)$ is exact at $H^iF$ if and only if $\beta_{i-1}=\ell_{i-1}$ and $\alpha_i=0$, i.e.  $\mathrm{can}^{i-1}$ is surjective and $\mathrm{can}^i$ is injective, which shows the equivalence of 2. and 3. using the first claim.
\end{proof}

\begin{rmk}
	Criterions 2. in the first claim and 3. in the second claim do not depend on the $B$-module structure, hence they can be checked on the underlying abelian groups or complexes thereof. 
\end{rmk}

\begin{cor}\label{cor:semisimplicity_frobenius}
	Assume $\mathbf{L}_{F,n}$ or $\mathbf{R}(\dim X)$.
	\begin{enumerate}
		\item The endomorphism $1-\phi$ of $R\Gamma_{\mathrm{c}}(X_{\overline{k}},F_\ell(n))_{\QQ_\ell}$ is semi-simple at $0$.
	\end{enumerate}
	We can thus write $R\Gamma_{\mathrm{c}}(X_{\overline{k}},F_\ell)_{\QQ_\ell}=V\oplus W$ where $\phi=1$ on $V$ and $1-\phi$ is an automorphism on $W$, inducing a decomposition  $R\Gamma_{\mathrm{c}}(X,F_\ell(n))_{\QQ_\ell}\simeq V[-1]\oplus V$.
	\begin{enumerate}[resume]
		\item The trivialization $\alpha_\ell$ is also given by
		\[
		[R\Gamma_{\mathrm{c}}(X,F_\ell(n))_{\QQ_\ell}]\to [V[-1]\oplus V] \to  [V]^{-1}[V] \xrightarrow{\mathrm{can}} \mathbf{1}. 
		\]
	\end{enumerate}
\end{cor}

\begin{proof}
	We know that $\cup e$ makes $H^\ast_{\mathrm{W,c}}(X,F(n))_\QQ$ into an acyclic complex, hence this is also true after base change to $\QQ_\ell$. The element $e_\ell:=e\otimes 1$ classifies the map $\beta\otimes 1:\ZZ_\ell \to \ZZ_\ell[1]$. It remains to identify the map $\beta_0$ of \cref{lem:semisimplicity} with $-(\beta\otimes 1)$; this is \cite[Lem. 1.2]{Rapoport1982}.\footnote{The sign comes from the choice of the arithmetic Frobenius to define $\beta$, while we are looking at the action the geometric Frobenius.}
	Then on the one hand we have $\beta_0=-(\beta\otimes 1)$ so that $(H^\ast_{\mathrm{c}}(X_{\overline{k}},F_\ell)_{\QQ_\ell},\cup (e\otimes 1))\xrightarrow[\simeq]{(-1)^\ast} (H^\ast_{\mathrm{c}}(X_{\overline{k}},F_\ell)_{\QQ_\ell},H^\ast \beta_0)$, and on the other hand under the decomposition $R\Gamma_{\mathrm{c}}(X,F_\ell(n))_{\QQ_\ell}\simeq V[-1]\oplus V$ we can identify $\beta_0$ as the composite map
	\[
	 V[-1]\oplus V \xrightarrow{\begin{pmatrix}	0 & 1\\ 0 & 0 \end{pmatrix}} V\oplus W \xrightarrow{\begin{pmatrix}	1 & 0\\ 0 & 0 \end{pmatrix}} V\oplus V[1]
	\]
	which is $\begin{pmatrix}	0 & 1 \\ 0 & 0	\end{pmatrix}.$	Repeating the argument of \cref{two_trivializations}, we obtain the second claim.
\end{proof}

\begin{rmk}
	Fix a prime $\ell\neq p$ and choose a compatible system of roots of unity in $k^{\mathrm{sep}}$. We obtain an isomorphism $\ZZ_\ell(n)\simeq \ZZ_\ell$ of étale sheaves over $\Spec(k^{\mathrm{sep}})$ under which the geometric Frobenius on the left hand side corresponds to multiplication by $q^{-n}$ on the right hand side.\footnote{This is a unit in $\ZZ_\ell$.} Thus for any $R_\ell$-constructible sheaf $G$, by the projection formula there is an induced identification $R\Gamma_{\mathrm{c}}(X_{\overline{k}},G(n))\simeq R\Gamma_{\mathrm{c}}(X_{\overline{k}},G)$ where the geometric Frobenius on the left corresponds to $\phi\cdot q^{-n}$ on the right. Therefore the corollary can be restated as the semi-simplicity of $1-\phi\cdot q^{-n}$ on $R\Gamma_{\mathrm{c}}(X_{\overline{k}},F_\ell)_{\QQ_\ell}$.
\end{rmk}

\begin{prop}
	If $F$ is a locally constant $R$-constructible complex and $X$ is smooth projective, there is a natural identification $V\simeq R\Gamma(X,F(n))_{\QQ_\ell}$. If $F$ is any $R$-constructible complex on a finite type separated scheme $X$ over $k$ and we assume $\mathbf{R}(\dim X)$, there is a natural identification $V\simeq R\Gamma_{\mathrm{eh},c}(X,F(n))_{\QQ_\ell}$. 
\end{prop}

\begin{proof}
	 We have an isomorphism $V\simeq \oplus H^i_{\mathrm{c}}(X_{\overline{k}},F_\ell(n))^{\phi=1}[-i]$, and there is thus a canonical $A_\ell$-linear map $R\Gamma(X,F(n))_{\QQ_\ell} \to V$ when $X$ is proper, resp. $R\Gamma_{\mathrm{eh,c}}(X,F(n))_{\QQ_\ell} \to V$ in the general case, induced through the semi-simplicity of $A_\ell$ by 
	 \[H^i(X,F(n))_{\QQ_\ell}\to H^i_{\mathrm{W}}(X,F(n))_{\QQ_\ell} \simeq H^i(X,F_\ell(n))_{\QQ_{\ell}}\to H^i(X_{\overline{k}},F_\ell(n))^{\phi=1},\]
	 resp.
	 \[
	 H^i_{\mathrm{eh,c}}(X,F(n))_{\QQ_\ell}\to H^i_{\mathrm{Wh,c}}(X,F(n))_{\QQ_\ell} \simeq H^i_{\mathrm{c}}(X,F_\ell(n))_{\QQ_{\ell}}\to H^i_{\mathrm{c}}(X_{\overline{k}},F_\ell(n))^{\phi=1}.
	 \]
	 Hence we can assume without loss of generality that $R=\ZZ$. Since the result is trivial if $F$ is a torsion constructible sheaf, we can apply \cref{lem:devissage} and, for the second claim, the étale--eh comparisons \cite[Thm. 3.6, Thm. 4.3]{Geisser2006} to reduce to the case where $F=\ZZ$ and $X$ is smooth projective.
	 
	 For $\ell=p$ and $n<0$, there is nothing to prove as everything is $0$ by definition. If $\ell\neq p$ and $n<0$, by Deligne's proof of the Weil conjectures \cite{Deligne1974} we find that $H^i(X_{\overline{k}},\QQ_\ell(n))$ is pure of weight $i-2n$, hence $1-\phi$ is an automorphism of $R\Gamma(X_{\overline{k}},\QQ_\ell(n))$ so that $V=0$; since $\ZZ(n)$ is torsion, we also have $R\Gamma(X,\ZZ(n))_{\QQ_\ell}=0$.
	 
	 For $n=0$, we obtain that $V=H^0(X_{\overline{k}},\QQ_\ell)^{\phi=1}[0]\simeq H^0(X,\QQ_\ell)[0]$. 
	We have $H^0(X,\ZZ)=\ZZ^{\#(\pi_0 X)}$, and for any normal Noetherian scheme we have $H^1(X,\ZZ)=\Hom_{\mathrm{cont}}(\pi_1^{\mathrm{proet}}(X),\ZZ)=0$ since $\pi_1^{\mathrm{proet}}(X)\simeq \pi_1^{\mathrm{et}}(X)$ is profinite, and $H^i(X,\QQ)=0$ for $i>0$ \cite[II, lem. 2.10]{Milne2006}. We thus obtain $\faktor{H^0(X,\ZZ)}{\ell^n} \simeq  H^0(X,\ZZ/\ell^n)$ and an isomorphism
	 \begin{align*}
	 R\Gamma(X,\ZZ)_{\QQ_\ell}\simeq H^0(X,\ZZ)_{\QQ_\ell}[0]\simeq (H^0(X,\ZZ)_{\ZZ_\ell})_{\QQ_\ell}[0]\simeq (H^0(X,\ZZ)^\wedge_\ell)_{\QQ_\ell}[0] &\simeq (\lim H^0(X,\ZZ/\ell^n))_{\QQ_\ell} [0]\\
	 &\simeq H^0(X,\QQ_\ell)[0]\\
	 &=H^0(X_{\overline{k}},\QQ_\ell)^{\phi=1}[0]\\
	 &=V.\qedhere
	 \end{align*}
\end{proof}

\section{Special values of cohomological \texorpdfstring{$L$}{L}-functions}\label{sec:special_value_cohomological}

In this section, we compute special values for characteristic polynomials or characteristic rational functions of semi-simple endomorphisms, and apply the resulting formula to the case of the characteristic rational function of the geometric Frobenius acting on the compactly supported $\ell$-adic cohomology of some coefficients. The main tool is \cref{lem:burns}, which is due to Burns and for which we provide a different proof.

\paragraph{} \label{par:local_case} Let $\ell$ be a prime, let $K$ be a finite extension of $\QQ_\ell$ with ring of integers $\cal{O}_K$, and let $B$ be an $\cal{O}_K$-order in a semi-simple algebra $B_K$ over $K$. We define the (inverse) characteristic rational fraction of a perfect complex with an endomorphism:

\begin{defi}
	Let $M\in \Perf(B)$ with an endomorphism $\theta$. The inverse characteristic rational function of $(M,\theta)$ over $B_K$ is
	\[
	\xi(\theta,T):=\mathrm{Nrd}_{B_K}(1-\theta\cdot T\mid M\otimes_B L_0 B_K[T])^{-1}\in (L_0\ct{B_K}[T])^\times
	\] 
	where $L_0(-)[T]$ is the Cohn localization of the polynomial ring in one variable at matrices which become invertible when evaluated at $0$.
\end{defi}

\begin{thm}[Special values at $\ell$-adic integers of characteristic rational functions]\label{thm:cohomological_special_value}
	Let $M\in \Perf(B)$ with an endomorphism $\theta$. Let $x\in \cal{O}_K$, and assume that $1-\theta\cdot x$ is semi-simple at $0$ on $M_K:=M\otimes_{\cal{O}_K}K$ with a $\theta$-equivariant decomposition $M_K=V\oplus W$ where $\theta\cdot x=1$ on $V$ and $1-\theta\cdot x$ is an automorphism on $W$. Then $\fib\left(M \xrightarrow{1-\theta\cdot x}M\right)\in \Perf(B)$.
	Define a trivialization of $[\fib(1-\theta\cdot x)_{K}]_{B_K}$ by 
	\[
	\alpha:[\fib(1-\theta\cdot x)_{K}]_{B_K}\xrightarrow{\simeq} [V[-1]\oplus V]_{B_K}\xrightarrow{\simeq} [V]_{B_K}^{-1}[V]_{B_K}\xrightarrow{\simeq} \mathbf{1}_{B_K}
	\]
	Then we have
	\begin{gather*}
		\ord_{t=x} \xi(\theta,T)=\sum_i (-1)^i\cdot i \cdot \mathrm{rrank}_{B_K} H^i(\fib(1-\theta\cdot x)_{K}) \\
		\delta(\xi^\ast(\theta,x))=-\chi([\fib(1-\theta\cdot x)]_{B},\alpha)
	\end{gather*}
	where $\delta:\ct{B_K}^\times\simeq K_1(B_K) \to K_0(B,B_K)$ is the boundary. Moreover, if $v_\ell(x)>0$, then $\ord_{t=x} \xi(\theta,T)=0$ and $\delta_\ell(\xi^\ast(\theta,x))=0$.
\end{thm}

\begin{rmk}
	\begin{itemize}
		\item []
		\item Let $N$ be the monoid generated by the geometric Frobenius inside the Galois group $\Gal(k^{\mathrm{sep}}/k)$. Define $\cal{O}_K\{x\}$ to be the $N$-equivariant étale sheaf on $\Spec(k^{\mathrm{sep}})$ corresponding to $\cal{O}_K$ with the monoid action of $N$ where $\phi$ acts by multiplication by $x$, and for any $B$-constructible sheaf $G$ on $X$, define an $N$-equivariant sheaf $G\{x\}$ on $X_{\overline{k}}$ by tensoring, over $\cal{O}_K$, the pullback of $G$ with the pullback to $X_{\overline{k}}$ of $\cal{O}_K\{x\}$. In the above situation, taking $M=R\Gamma_{\mathrm{c}}(X_{\overline{k}},G)$ and $\theta=\phi$, the above complex $\fib(1-\phi\cdot x)$ then computes the natural cohomology with compact support of $G\{x\}$ on $X$.
		\item If $\cal{O}_K$ contains a root $x$ of $X^2=q$, this gives a formula for the  special value at $x=q^{1/2}$, and also of the special value at $x^{-1}=q^{-1/2}$ when $\ell\neq p$.
		\item The argument works more generally for $B$ an order over a complete discrete valuation ring as soon as the reduced norm $\mathrm{Nrd}_{B_K}:K_1(B_K)\to \ct{B_K}^\times$ is an isomorphism.
	\end{itemize}
\end{rmk}

\begin{proof}
	From the $\theta$-equivariant decomposition $M_{K}=V\oplus W$
	we get
	\begin{align*}
		\xi(\theta,T)&=\mathrm{Nrd}_{B_K}(1-x^{-1}\cdot T\mid V\otimes_{B_K}B_K[[T]])^{-1}\cdot\mathrm{Nrd}_{B_K}(1-\theta\cdot T\mid W\otimes_{B_K}B_K[[T]])^{-1}\\
		&=x^{-\mathrm{rrank}_{B_K} V}\cdot\mathrm{Nrd}_{B_K}(x- T\mid V\otimes_{B_K}B_K[[T]])^{-1}\cdot\mathrm{Nrd}_{B_K}(1-\theta\cdot T\mid W\otimes_{B_K}B_K[[T]])^{-1}.
	\end{align*}
	For an associative algebra $B$ over $\cal{O}_K$, denote by $L_{0,x} B[T]$ the Cohn localization of $B$ at matrices that become invertible when evaluated at $0$ and $x$. For the commutative ring $\ct{B_K}[T]$, $L_{0,x} \ct{B_K}[T]$ is the ring of rational functions without pole at $0$ or $x$ (on each component). Since the endomorphism $1-\theta\cdot T$ of $W\otimes_{B_K} B_K[T]$ is an automorphism when evaluated at $0$ and $x$, it defines a class 
	\[
	\langle 1-\theta\cdot T\mid W\otimes_{B_K} L_{0,x} B_K[T]\rangle \in K_1(L_{0,x} B_K[T])
	\] whose reduced norm, which is naturally an element of $\left(L_{0,x}(\ct{B_K}[T])\right)^\times$, coincides with $\mathrm{Nrd}_{B_K}(1-\theta\cdot T\mid W\otimes_{B_K} B_K[[T]])$ inside $\ct{B_K}[[T]]^\times$; in other words, $\mathrm{Nrd}_{B_K}(1-\theta\cdot T\mid W\otimes_{B_K} B_K[[T]])$ is rational without zeros or poles at $0$ and $x$. It follows that the order of vanishing is determined by $\mathrm{Nrd}_{B_K}(1-T\mid V\otimes_{B_K} B_K[[T]])^{-1}$ and that the special value at $x$ is given by 
	\[
	\xi^\ast(\theta,x)=x^{-\mathrm{rrank}_{B_K} V}\cdot\mathrm{ev}_{t=x}(\mathrm{Nrd}_{B_K}( 1-\theta\cdot T\mid W\otimes_{B_K} L_{0,x} B_K[T])^{-1}).
	\]
	The reduced norm is, after base change to $\overline{K}$, a determinant, so it commutes with evaluation at $t=x$, hence:
	\begin{align*}
		\mathrm{ev}_{t=x}(\mathrm{Nrd}_{B_K}( 1-\theta\cdot T\mid W\otimes_{B_K} L_{0,x} B_K[T])^{-1})&=\mathrm{Nrd}_{B_K}(\mathrm{ev}_{t=x}(\langle 1-\theta\cdot T\mid W\otimes_{B_K} L_{0,x} B_K[T]\rangle))^{-1}\\
		&=\mathrm{Nrd}_{B_K}(1-\theta\cdot x\mid W)^{-1},
	\end{align*}
	which identifies with $\langle 1-\theta\cdot x\mid W \rangle^{-1}\in K_1(B_K)$ under the isomorphism $K_1(B_K)\simeq\ct{B_K}^\times$. We obtain:
	\[
	\delta(\xi^\ast(\theta,x))=-\delta(\mathrm{Nrd}_{B_K}(1-\theta\cdot x\mid W))+\delta(x^{-\mathrm{rrank}_{B_K} V}).
	\]
	
	Now we observe that
	\begin{align*}
		\ord_{t=x}\mathrm{Nrd}_{B_K}(x-T\mid V)^{-1}=-\mathrm{rrank}_{B_K}V &=-\sum_{i\in \ZZ} (-1)^i\mathrm{rrank}_{B_K} H^iV\\
		& = \sum_{i\in \ZZ} (-1)^i\cdot i \cdot \mathrm{rrank}_{B_K} H^i(\fib(1-\theta\cdot x)_{K})
	\end{align*}
	by telescoping the sum, since $ H^i(\fib(1-\theta\cdot x)_{K})=H^{i-1}V\oplus H^iV$.
	
	If $v_\ell(x)>0$, then $1-\theta\cdot x$ is an automorphism of $M$, so that we must have $V=0$ and thus $\ord_{t=x}\xi(\theta,T)=-\mathrm{rrank}_{B_K}V=0$. This implies that $x^{-\mathrm{rrank}_{B_K}V}=1$ so that the special value is then just $\langle 1-\theta\cdot x\mid W \rangle^{-1}$, which lifts to $K_1(B)$ and thus vanishes under $\delta$. On the other hand, if $v_\ell(x)=0$, then $x\in \ct{B_K}^\times$ lifts to $K_1(B)$ and thus vanishes under $\delta$. We find in particular in both cases that $\delta(x^{-\mathrm{rrank}_{B_K}V})=0$.
	
	It remains to show that $\delta(\langle 1-\theta\cdot x\mid W \rangle^{-1})=-\chi([\fib(1-\theta\cdot x)]_{B},\alpha)$. This will follow through the next lemma.
\end{proof}

\begin{lem}[{\cite[Prop. 3.2 (ii)]{Burns}}]\label{lem:burns}
	Let $B$ be an associative ring with unit, $B\to C$ a ring morphism, and let $M\in \Perf(B)$ with an endomorphism $\theta$. Denote by $F$ the fiber of $\theta$. Assume that $\theta$ is semi-simple at $0$ on $M\otimes_B C$, with a $\theta$-equivariant decomposition $M\otimes_B C\simeq V\oplus W$ where $\theta=0$ on $V$ and $\theta$ is an automorphism of $W$. Let
	\[
	\alpha:[F\otimes_B C]_C \to [V[-1]\oplus V]_C \to [V]_C^{-1}[V]_C\to \mathbf{1}_C
	\]
	be the trivialization of $[F\otimes_B C]_C$ induced by the direct sum decomposition coming from the semi-simplicity. Denote by $\delta$ the co-boundary map $K_1(C)\to K_0(B,C)$. Then we have an equality
	\[
	\chi([F]_B,\alpha)=\delta(\langle \theta \mid W\rangle)\in  K_0(B,C).
	\]
\end{lem}

\begin{proof}
	Recall that $\delta$ is modeled by the map sending $\gamma\in \mathrm{Aut}(\mathbf{1}_C)$ to the element $(\mathrm{1}_B,\gamma)$ of the relative Picard groupoid $V(B,C)$, and that if $\psi$ is an automorphism of a perfect complex $Y\in \Perf(C)$, then $\langle \psi \mid Y\rangle \in K_1(C)$ is the element $\mathrm{can}\circ ([\psi]_C\otimes 1)\circ \mathrm{can}^-1\in \mathrm{Aut}(\mathbf{1}_C)$, where $\mathrm{can}:[Y]_C[Y]_C^{-1}\to \mathbf{1}_C$ is the map witnessing $[Y]_C^{-1}$ as the inverse of $[Y]_C$.
	
	Since $M\in \Perf(B)$, there is a canonical trivialization of $F$ given by $\varepsilon:[F]_B \to [M]_B[M]_B^{-1}\to \mathbf{1}_B$. We thus get an isomorphism $\varepsilon:([F],\alpha)\xrightarrow{\simeq} (\mathbf{1}_B,\alpha\circ\varepsilon_C^{-1})$ in the relative Picard groupoid $V(B,C)$, hence $\chi([F]_B,\alpha)=\chi(\mathbf{1}_B,\alpha\circ\varepsilon_C^{-1})=\delta(\alpha\circ\varepsilon_C^{-1})$. We now compute $\alpha\circ\varepsilon_C^{-1}$.
	Let $a$ and $b$ be the following fiber sequences:
	\[
	\begin{array}{llc}
		(a) & & V[-1]\oplus V \to V \xrightarrow{0} V\\
		(b) & & 0 \to W \xrightarrow{\theta} W
	\end{array}
	\]
	We also have an isomorphism of fiber sequences
	\[\begin{tikzcd}[ampersand replacement=\&]
		{F\otimes_B C} \& {M\otimes_B C} \& {M\otimes_B C} \\
		{V[-1]\oplus V} \& {V\oplus W} \& {V\oplus W}
		\arrow[from=1-1, to=1-2]
		\arrow["\simeq", from=1-1, to=2-1]
		\arrow["{\theta\otimes 1}", from=1-2, to=1-3]
		\arrow["\simeq", from=1-2, to=2-2]
		\arrow["\simeq", from=1-3, to=2-3]
		\arrow[from=2-1, to=2-2]
		\arrow["{\begin{pmatrix}0&0 \\ 0 & \theta_{\mid W}\end{pmatrix}}", from=2-2, to=2-3]
	\end{tikzcd}\]
	Consider the following commutative diagram in the Picard groupoid $V_C$:
	\[\begin{tikzcd}[ampersand replacement=\&]
		{[F\otimes_B C]} \& {[M\otimes_B C][M\otimes_B C]^{-1}} \& {\mathbf{1}_C} \\
		{[V[-1]\oplus V]} \& {[V\oplus W][V\oplus W]^{-1}} \& {\mathbf{1}_C} \\
		{[V[-1]\oplus V]} \& {[V][V]^{-1}[W][W]^{-1}} \& {\mathbf{1}_C} \\
		{[V][V]^{-1}} \& {[V][V]^{-1}[W][W]^{-1}} \& {\mathbf{1}_C} \\
		{\mathbf{1}_C} \& {[W][W]^{-1}} \& {\mathbf{1}_C}
		\arrow[from=1-1, to=1-2]
		\arrow[from=1-1, to=2-1]
		\arrow["{\mathrm{can}}", from=1-2, to=1-3]
		\arrow[from=1-2, to=2-2]
		\arrow[Rightarrow, no head, from=1-3, to=2-3]
		\arrow[from=2-1, to=2-2]
		\arrow[Rightarrow, no head, from=2-1, to=3-1]
		\arrow["{\mathrm{can}}"', from=2-2, to=2-3]
		\arrow[from=2-2, to=3-2]
		\arrow[Rightarrow, no head, from=2-3, to=3-3]
		\arrow["{[a]\otimes[b]}"', from=3-1, to=3-2]
		\arrow["{[a]}"', from=3-1, to=4-1]
		\arrow["{\mathrm{can}\otimes \mathrm{can}}"', from=3-2, to=3-3]
		\arrow[Rightarrow, no head, from=3-2, to=4-2]
		\arrow[Rightarrow, no head, from=3-3, to=4-3]
		\arrow["{1\otimes[b]}"', from=4-1, to=4-2]
		\arrow["{\mathrm{can}}"', from=4-1, to=5-1]
		\arrow["{\mathrm{can}\otimes \mathrm{can}}"', from=4-2, to=4-3]
		\arrow["{\mathrm{can}\otimes 1}", from=4-2, to=5-2]
		\arrow[Rightarrow, no head, from=4-3, to=5-3]
		\arrow["{[b]}"', from=5-1, to=5-2]
		\arrow["{\mathrm{can}}"', from=5-2, to=5-3]
	\end{tikzcd}\]
	where the upper part is induced from the above isomorphism of fiber sequences. The left vertical composite is $\alpha$, the top horizontal composite is $\varepsilon_C$, hence $\alpha\circ\varepsilon_C^{-1}=[b]^{-1}\mathrm{can}^{-1}$. The following isomorphism of fiber sequences
	\[\begin{tikzcd}[ampersand replacement=\&]
		{0 } \& W \& W \\
		0 \& W \& W
		\arrow[from=1-1, to=1-2]
		\arrow[equal, from=1-1, to=2-1]
		\arrow["\theta_{\mid W}"', from=1-2, to=1-3]
		\arrow["\theta_{\mid W}", from=1-2, to=2-2]
		\arrow[equal, from=1-3, to=2-3]
		\arrow[from=2-1, to=2-2]
		\arrow[equal, from=2-2, to=2-3]
	\end{tikzcd}\]
	induces a commutative diagram
	\[\begin{tikzcd}[ampersand replacement=\&]
		{\mathbf{1}_C} \& {[W][W]^{-1}} \\
		{\mathbf{1}_C} \& {[W][W]^{-1}}
		\arrow["{[b]}", from=1-1, to=1-2]
		\arrow[Rightarrow, no head, from=1-1, to=2-1]
		\arrow["{[\theta]\otimes 1}", from=1-2, to=2-2]
		\arrow["{\mathrm{can}^{-1}}"', from=2-1, to=2-2]
	\end{tikzcd}\]
	so that 
	\[
	\alpha\circ\varepsilon_C^{-1}=[b]^{-1}\mathrm{can}^{-1}=\mathrm{can}\circ\mathrm{can}^{-1}\circ[b]^{-1}\circ\mathrm{can}^{-1}=\mathrm{can}\circ([\theta_{\mid W}]\otimes 1)\circ\mathrm{can}^{-1}=\langle \theta \mid W\rangle
	\]
\end{proof}

\paragraph{} We now return to characteristic rational fractions associated to constructible sheaves. We keep notations from the previous paragraph. Let $G$ be a $B$-constructible complex on a finite type separated scheme $X$ over $k$. We first treat the case where the residual characteristic of $\cal{O}_K$ is prime to $p$.

\begin{cor}[Special values in terms of Tate twists for $\ell$-adic sheaves]\label{cor:special_value_ell}
	Let $G$ be a $B$-constructible complex on $X$ such that $\cal{O}_K$ has residual characteristic $\ell\neq p$, let $n\in \ZZ$ and assume that $1-\phi$ is semi-simple at $0$ on $R\Gamma_{\mathrm{c}}(X_{\overline{k}},G(n))_{K}$. We have
	\begin{gather*}
		\ord_{t=q^{-n}} \xi_X(G,T)=\sum_i (-1)^i\cdot i \cdot \mathrm{rrank}_{B_K} H^i_{\mathrm{c}}(X,G(n))_{K} 	\\
		\delta(\xi_X^\ast(G,q^{-n}))=-\chi_X(G(n))
	\end{gather*}
\end{cor}

\begin{proof}
	  We have for the pair $(R\Gamma_{\mathrm{c}}(X_{\overline{k}},G),\phi)$ that $\xi(\phi,T)=\xi_X(G,T)$. We have observed that there is an identification $R\Gamma_{\mathrm{c}}(X_{\overline{k}},G(n))\simeq R\Gamma_{\mathrm{c}}(X_{\overline{k}},G)$ where the geometric Frobenius on the left corresponds to $\phi\cdot q^{-n}$ on the right. Under this identification, we have 
	 \[
	 R\Gamma_{\mathrm{c}}(X,G(n))\simeq \fib\left(R\Gamma_{\mathrm{c}}(X,G)\xrightarrow{1-\phi\cdot q^{-n}}R\Gamma_{\mathrm{c}}(X,G)\right)
	 \]
	 hence everything follows from \cref{thm:cohomological_special_value}.
\end{proof}

We now treat the case where the residual characteristic of $\cal{O}_K$ is $p$; the same proof as above yields:
\begin{cor}[Special values at $1$ for $p$-adic sheaves]\label{cor:special_value_p_1}
	Let $G$ be a $B$-constructible complex on $X$ such that $\cal{O}_K$ has residual characteristic $p$ and assume that $1-\phi$ is semi-simple at $0$ on $R\Gamma_{\mathrm{c}}(X_{\overline{k}},G)_{K}$. We have
	\begin{gather*}
		\ord_{t=1} \xi_X(G,T)=\sum_i (-1)^i\cdot i \cdot \mathrm{rrank}_{B_K} H^i_{\mathrm{c}}(X,G)_{K} 	\\
		\delta(\xi_X^\ast(G,1))=-\chi_X(G)
	\end{gather*}
\end{cor}

Tate twists do not make sense for $p$-adic sheaves in characteristic $p$. However, when $n<0$, $1-\phi\cdot q^{-n}$ is invertible on $R\Gamma_c(X_{\overline{k}},G)$ by $p$-completeness, so that 
\[\fib(R\Gamma_c(X_{\overline{k}},G)\xrightarrow{1-\phi\cdot q^{-n}} R\Gamma_c(X_{\overline{k}},G))=0.
\]
Comparing this with the $\ell$-adic case where the output is $R\Gamma_{\mathrm{c}}(X,G(n))$, this justifies to define $G(n)=0$, which is compatible with our previous definition when $G$ is the $p$-adic completion of an integrally defined constructible sheaf. Moreover, this is compatible with the expected special value formula:
\begin{cor}[Special values at negative integers\footnote{Usually, special values are investigated for $L$-functions in the $s$-variable; here $T=q^{-s}$ so negative integers correspond to positive powers of $q$.} for $p$-adic sheaves]\label{cor:special_value_p_negative}
	Let $n<0$. We have
		\begin{gather*}
		\ord_{t=q^{-n}} \xi_X(G,T)=0=\sum_i (-1)^i\cdot i \cdot \mathrm{rrank}_{B_K} H^i_{\mathrm{c}}(X,G(n))_{K} 	\\
		\delta(\xi_X^\ast(G,q^{-n}))=0=-\chi_X(G(n))
	\end{gather*}
\end{cor}
\begin{proof}
	We can apply \cref{thm:cohomological_special_value} to the pair $(R\Gamma_{\mathrm{c}}(X_{\overline{k}},G),\phi)$ since the remark above implies the semisimplicity at $0$ of $1-\phi\cdot q^{-n}$ on $R\Gamma_{\mathrm{c}}(X_{\overline{k}},G)$. It remains to check that the right hand side in the two formulas vanish, but this follows from $G(n)$ being $0$.
\end{proof}

\section{Proof of the special value theorem}\label{sec:proof}

\begin{thm}\label{thm:main_body}
	Let $X$ be a finite type separated scheme over a finite field $k=\FF_q$, let $R$ be an order in a semi-simple $\QQ$-algebra $A$, and let $F$ be an $R$-constructible complex on $X$. Let $n\leq 0$.
	\begin{itemize}
		\item Assume $\mathbf{L}_{F,n}$. Then $R\Gamma_{\mathrm{W,c}}(X,F(n))\in\Perf(R)$, we have the vanishing order formula\footnote{If $n<0$, $F(n)$ is torsion so the formula says that the vanishing order is zero}
		\[
		\mathrm{ord}_{t=q^{-n}}L_X(F,T)=\sum_i (-1)^i\cdot i \cdot \mathrm{rrank}_A H^i_{\mathrm{W,c}}(X,F(n))_\QQ
		\]
		and the special value formula
		\[
		\widehat{\delta_R} (L^\ast_X(F,q^{-n}))=-\chi_{\mathrm{W},X}(F(n)).
		\]
		\item Assume $\mathbf{R}(\dim X)$. The similar result holds using $R\Gamma_{\mathrm{Wh,c}}(X,F(n))$ .
	\end{itemize}
\end{thm}
\begin{proof}
Let us treat the claim for Weil-étale cohomology. Thus, assume $\mathbf{L}_{F,n}$. Recall that $R\Gamma_{\mathrm{W,c}}(X,F(n))\in\Perf(R)$ by \cref{cor:perfect_Z_implies_perfect_R} whence the Weil-étale Euler characteristic $\chi_{\mathrm{W},X}(F(n))$ is well-defined. By \cref{cor:reduction_local_case}, we have to check that
\[
\ord_{t=q^{-n}} \xi_X(F_\ell,T)= \sum_i (-1)^i\cdot i \cdot \mathrm{rrank}_{A_\ell} H^i_{\mathrm{c}}(X,F_\ell(n))_{\QQ_\ell} 
\]
for some $\ell\neq p$ and
\[
\delta_\ell(\xi_X^\ast(F_\ell,q^{-n}))=-\chi_X(F_\ell(n)),
\]
for \emph{all} primes $\ell$, where $F_\ell$ is the $\ell$-adic completion of $F$ and $\delta_\ell$ is the boundary $\ct{A_\ell}^\times\simeq K_1(A_\ell) \to K_0(R_\ell,A_\ell)$.
Those formulas then hold by the computation of special values of cohomological $L$-functions \cref{cor:special_value_ell,cor:special_value_p_1,cor:special_value_p_negative} which can be applied because Weil-étale cohomology is an integral model of $\ell$-adic cohomology and thus forces semi-simplicity at $0$ of $1-\phi$ \cref{cor:semisimplicity_frobenius}. The latter corollary also implies that the two trivializations of \cref{def:local_Euler_char} and \cref{thm:cohomological_special_value} agree, so they define the same refined Euler characteristic.
\end{proof}

\section{Special values for non-commutative \texorpdfstring{$L$}{L}-functions of motives}\label{sec:applications}

In this final section, we apply our results in the case where we can identify the non-commutative $L$-function of a (derived) motive over a global field $K$ of characteristic $p$, with action of a semi-simple $\QQ$-algebra, with the $L$-function of a suitable constructible complex on a finite type separated scheme over $\FF_p$. This can be done for Artin motives (\cref{prop:constructible_model_Artin_motive}) or for total derived motives of smooth projective varieties over $K$ that have a smooth projective model over $\FF_p$, leading to special value formulas for those $L$-functions (\cref{cor:special_value_artin,thm:special_value_motive}).

Let $K$ be a global function field, let $S$ be a finite set of places of $K$ and let $V$ be a finite $A$-module with a discrete action of $G_K$, which we identify with an $A$-linear Artin motive over $K$. Denote by $C$ the smooth projective curve with function field $K$ and $U:=C\backslash S$ its open subscheme. Denote by $g:\Spec(K)\to U$ the inclusion of the generic point. 

\begin{prop}\label{prop:constructible_model_Artin_motive}
	We have:
	\begin{enumerate}
		\item $g_\ast V$ is an $A$-constructible sheaf on $U$;
		\item There is a flat $R$-constructible sheaf $F$ on $U$ such that $F\otimes \QQ=g_\ast V$.
	\end{enumerate} 
\end{prop}

\begin{proof}
	The sheaf $g_\ast V$ is locally constant on a dense open, with stalk at the generic point the finite $A$-module $V$. At the finitely many remaining closed points $x$ the stalk is $(g_\ast V)_x\simeq V^{I_x}$ which is also finite, so $g_\ast V$ is $A$-constructible by semisimplicity.
	
	The action of $G_K$ on $V$ factors through a finite quotient $G$, and $V$ is a finite $A[G]$-module. By semisimplicity, $V$ is a direct factor of some $A[G]^n$, so is determined by an idempotent $e\in M_n(A[G])$. By clearing denominators we can assume that $e\in M_n(R[G])$, and thus $e$ determines a finite projective $R[G]$-module $M$ with $M\otimes\QQ=V$. In particular, $M$ is a finite projective $R$-module with a discrete action of $G_K$. Then $g_\ast M$ is locally constant $R$-constructible on a dense open in $U$; it remains to correct $g_\ast M$ at the finitely many points where the ramification acts non-trivially on $V$ so that we get an $R$-constructible sheaf $F$ on $U$. Let $T$ denote the set of closed points $x$ where the inertia group $I_x$ acts non-trivially on $V$, and let $G_x$ denote the Galois group of the residue field at $x$. It suffices to produce, for each $x\in T$, a finite projective $R$-module $M_x$ with an action of $G_x$ and a $G_x$-equivariant map $M_x \to M^{I_x}$ that identifies with $V^{I_x}\xrightarrow{\mathrm{id}} V^{I_x}$ after tensoring with $\QQ$. We first find as before a finite projective $R$-module $M_x$ with a discrete action of $G_x$ such that $M_x\otimes \QQ\simeq V^{I_x}$. Now each generator of $M_x$ can be written as a $\QQ$-linear combination of generators of $M^{I_x}$ so that for some integer $d>0$ the composite $G_x$-equivariant map $M_x\hookrightarrow V^{I_x} \xrightarrow{d} V^{I_x}$ factors through the inclusion $M^{I_x} \hookrightarrow V^{I_x}$ as an $R$-linear $G_x$-equivariant map. We are done by replacing $M_x$ by its isomorphic image in $V^{I_x}$ under $M_x\hookrightarrow V^{I_x} \xrightarrow{d} V^{I_x}$.
\end{proof}

\begin{rmk}
	\begin{itemize}
		\item[]
		\item If $R$ is regular, we can take $F=g_\ast M$ for $M$ any finite projective $R$-module with a discrete action of $G_K$ such that $M\otimes \QQ=V$, but $g_\ast M$ will not be flat in general. 
		\item If $M$ is any finite projective $R$-module with a discrete action of $G_K$ such that $M\otimes \QQ=V$, and $M$  is unramified at any point of $U$, we can take $F=g_\ast M$.
	\end{itemize}
\end{rmk}

Since $\mathbf{L}_{F,n}$ holds for $n\leq 0$ and $F$ a $\ZZ$-constructible sheaf on a curve, we obtain:

\begin{cor}[Special values for non-commutative Artin $L$-functions]\label{cor:special_value_artin}
	Let $F$ be an $R$-constructible sheaf on $U$ such that $F\otimes\QQ=g_\ast V$. We have for the non-commutative Artin $L$-function $L_S(V,T)$ omitting the places in $S$:
	\[
	L_S(V,T)=L_U(F,T)
	\]
	and thus the vanishing order and special value formulas for any $n\leq 0$:
	\begin{gather*}
		\mathrm{ord}_{t=q^{-n}} L_S(V,T)= \sum_i (-1)^i\cdot i \cdot \mathrm{rrank}_{A} H^i_{\mathrm{W,c}}(U,F(n))_{\QQ}=\sum_i (-1)^i\cdot i \cdot \mathrm{rrank}_{A} H^i_{\mathrm{Wh,c}}(U,F(n))_{\QQ},\\
		\widehat{\delta_R}(L_S^\ast(V,q^{-n}))=-\chi_{W,U}(F(n))=-\chi_{Wh,U}(F(n)).
	\end{gather*}
\end{cor}

More generally, let $M$ be an $A$-linear object in the bounded derived category of pure motives over $K$. Let us formulate a technical hypothesis:

\begin{ass*}[$\mathbf{H}_M$]
	There is an $R$-constructible complex $F$ on a finite type separated scheme $f:X\to U$ such that 
	\[
	R^if_!(F_\ell[\frac 1 \ell])=g_\ast H^i\mathrm{Real}_\ell(M).
	\]
	for all $i\in \ZZ$, where $\mathrm{Real}_\ell$ is the derived $\ell$-adic realization, seen as a lisse $A_\ell$-constructible complex on $\Spec(K)$.
\end{ass*}

For instance, if $L/K$ is a Galois extension with group $G$, unramified over $U$ hence inducing a $G$-cover $V\to U$, then for $X_K$ a smooth projective variety over $K$ with smooth projective or more generally smooth proper model $f:X\to U$, the induced motive to $K$ of the total motive $M=\oplus h^i(X_L,0)[-i]$ of $X_L$ over $L$ is a $\QQ[G]$-linear derived motive over $K$; letting $\pi: X_V\to X$ denote the canonical projection, the sheaf $F:=\pi_\ast \ZZ$ with its natural $\ZZ[G]$-action satisfies the above assumption for $M$. Indeed, $F$ is a locally constant $\ZZ[G]$-constructible sheaf and $F_\ell[\frac 1 \ell]=\pi_\ast \QQ_\ell$, hence by \cite[Cor. VI.4.2]{Milne1980} $R^if_\ast \pi_\ast \QQ_\ell$ is a lisse $\QQ_\ell[G]$-sheaf. Since a lisse sheaf on $U$ is equal to its pull-push from the generic point, proper base change implies that
\[
R^if_\ast \pi_\ast \QQ_\ell=g_\ast g^\ast R^i(f\circ \pi)_\ast \QQ_\ell=g_\ast \oplus_G H^i(X_{\overline{K}},\QQ_\ell).
\]
More generally, if there is a smooth proper variety $X$ over $U$ and a locally constant $R$-constructible sheaf $F$ on $X$ such that the cohomology of the restriction of $F_\ell[\frac 1 \ell]$ to the generic fiber identifies with the $\ell$-adic realization of $M$, then $\mathbf{H}_M$ holds using the sheaf $F$. Note that this implies that $M$ has good reduction over $U$. The author does not see a systematic way to produce a sheaf $F$ outside of the good reduction hypothesis, except for Artin motives as investigated above.

\begin{thm}[Special values for motives of ``constructible nature'']\label{thm:special_value_motive}
	In the above setting, assume $\mathbf{H}_M$ and denote by $F$ the provided $R$-constructible sheaf. We have for the non-commutative motivic $L$-function $L_S(M,T)$ of $M$ avoiding the places in $S$:
	\[
	L_S(M,T)=L_X(F,T).
	\]
	Let $n\leq 0$, and
	\begin{enumerate}
		\item assume furthermore $\mathbf{L}_{F,n}$; then we have the vanishing order and special value formulas at $t=q^{-n}$:
			\begin{gather*}
				\mathrm{ord}_{t=q^{-n}} L_S(M,T)= \sum_i (-1)^i\cdot i \cdot \mathrm{rrank}_{A} H^i_{\mathrm{W,c}}(X,F(n))_{\QQ},\\
				\widehat{\delta_R}(L_S^\ast(M,q^{-n}))=-\chi_{W,X}(F(n)).
			\end{gather*}
			In particular, this holds unconditionally if the provided sheaf $F$ is locally constant and the associated scheme $X$ is smooth projective over $C$.
		\item assume furthermore $\mathbf{R}(\dim X)$; then we have the vanishing order and special value formulas at $t=q^{-n}$:
		\begin{gather*}
			\mathrm{ord}_{t=q^{-n}} L_S(M,T)= \sum_i (-1)^i\cdot i \cdot \mathrm{rrank}_{A} H^i_{\mathrm{Wh,c}}(X,F(n))_{\QQ},\\
			\widehat{\delta_R}(L_S^\ast(M,q^{-n}))=-\chi_{Wh,X}(F(n)).
		\end{gather*}
		In particular, this holds unconditionally if the provided sheaf $F$ is locally constant and the associated scheme $X$ is smooth proper over $C$ of dimension less than $2$, or less than $3$ for $p\geq 7$.
	\end{enumerate} 
\end{thm}

\printbibliography

\end{document}